%% file: paper.tex
\documentclass[final,1p,times]{elsarticle}

\usepackage{amssymb}
\usepackage{amsthm}
\usepackage{lineno}

\usepackage{amsmath}
\usepackage[ruled,algo2e,linesnumbered,algonl]{algorithm2e}
\newcommand{\vect}[1]{\boldsymbol{\mathbf{#1}}}
\usepackage[]{xcolor}

\definecolor{darkgreen}{rgb}{0.0, 0.35, 0.0}
\newcommand{\reva}[1]{\textcolor{black}{#1}}
\newcommand{\revb}[1]{\textcolor{black}{#1}}
\newcommand{\revc}[1]{\textcolor{black}{#1}}

\usepackage[free-standing-units,group-separator={,}]{siunitx}
\graphicspath{{figures/}}
\usepackage{subcaption}
\usepackage{url}

\journal{Journal of Computational Physics}

\begin{document}

\begin{frontmatter}



\title{An arbitrary order time-stepping algorithm for tracking particles in inhomogeneous magnetic fields}


\author[leeds]{Krasymyr Tretiak}
\ead{k.tretiak@leeds.ac.uk}
\author[leeds]{Daniel Ruprecht\corref{cor}}
\ead{d.ruprecht@leeds.ac.uk}

\address[leeds]{School of Mechanical Engineering, University of Leeds, United Kingdom}

\cortext[cor]{Corresponding author}

\begin{abstract}
The Lorentz equations describe the motion of electrically charged particles in electric and magnetic fields and are used widely in plasma physics.
The most popular numerical algorithm for solving them is the Boris method, a variant of the St\"ormer-Verlet algorithm.
Boris' method is phase space volume conserving and simulated particles typically remain near the correct trajectory.
However, it is only second order accurate.
Therefore, in scenarios where it is not enough to know that a particle stays on the right trajectory but one needs to know where on the trajectory the particle is at a given time, Boris method requires very small time steps to deliver accurate phase information, making it computationally expensive.
We derive an improved version of the high-order Boris spectral deferred correction algorithm (Boris-SDC) by adopting a convergence acceleration strategy for second order problems based on the Generalised Minimum Residual (GMRES) method.
Our new algorithm is easy to implement as it still relies on the standard Boris method.
Like Boris-SDC it can deliver arbitrary order of accuracy through simple changes of runtime parameter but possesses better long-term energy stability.
We demonstrate for two examples, a magnetic mirror trap and the Solev'ev equilibrium, that the new method can deliver better accuracy at lower computational cost compared to the standard Boris method.
While our examples are motivated by tracking ions in the magnetic field of a nuclear fusion reactor, the introduced algorithm can potentially deliver similar improvements in efficiency for other applications.
\end{abstract}

\begin{keyword}
Boris integrator\sep
particle tracking\sep
high-order time integration\sep
spectral deferred corrections\sep
fusion reactor



\end{keyword}

\end{frontmatter}


\input{intro.tex}

\input{method.tex}

%
%
\section{Numerical results}
We compare the performance of \revb{BGSDC} against both the staggered and non-staggered Boris method for two fusion-related test problems.
The first is a magnetic mirror where particles are confined by a magnetic field generated by two coils.
Our second benchmark uses a Solev'ev equilibrium which resembles the magnetic field in the Joint European Torus (JET) experimental Tokamak reactor.


\input{mtresults.tex}

\input{solresults.tex}

\input{conclusions.tex}

\section*{Acknowledgements}
This work was support by the Engineering and Physical Sciences Research Council EPSRC under grant EP/P02372X/1 ``A new algorithm to track fast ions in fusion reactors''. The authors thankfully acknowledge the very helpful input from Dr Rob Akers and James Buchanan from the Culham Centre for Fusion Energy and Dr Jonathan Smith from Tech-X UK Ltd.


\bibliographystyle{elsarticle-num} 
\bibliography{sdc,biblio,pint}


\end{document}

%% file: intro.tex
\section{Introduction}
The Lorentz equations
\begin{subequations}\label{eq:motion_equations}
	\begin{align}
	\vect{\dot{x}} &= \vect{v}, \\
	\vect{\dot{v}} &= \alpha \left[\vect{E}(\vect{x},t) + \vect{v} \times \vect{B}(\vect{x},t))\right] =: \vect{f}(\vect{x}, \vect{v})
	\end{align}
\end{subequations}
model movement of charged particles in electro-magnetic fields.
Here, $\vect{x}(t)$ is a vector containing all particles' position at some time $t$, $\vect{v}(t)$ contains all particles' velocities, $\alpha$ is the charge-to-mass ratio, $\vect{E}$ the electric field (both externally applied and internally generated from particle interaction) and $\vect{B}$ the magnetic field.
The Lorentz equations are used in many applications in computational plasma physics, for example laser-plasma interactions~\cite{Gibbon2005}, particle accelerators~\cite{VayEtAl2012} or nuclear fusion reactors~\cite{AkersEtAl2002}.

The Boris method, introduced by Boris in 1970~\cite{Boris1970}, is the most popular numerical scheme used for solving~\eqref{eq:motion_equations} although other numerical time stepping methods like Runge-Kutta-4 are used as well.
It is based on the Leapfrog algorithm but uses a special trick to resolve the seemingly implicit dependence that arises from the fact that the Lorentz force depends on $\vect{v}$.
Its popularity is because it is computationally cheap, second order accurate and phase space conserving~\cite{QinEtAl2013}.
While it was recently shown that for general magnetic fields this does not guarantee energy conservation and that the Boris method can exhibit energy drift~\cite{HairerLubich2018}, it is nevertheless a surprisingly good algorithm.
In most cases, particles will remain close to their correct trajectory because of its conservation properties.
However, Boris' method can introduce substantial phase errors and, for long time simulations, it only ensures that particles are near the right trajectory -- it does not provide information about where on the trajectory they are at a given time.

For some applications this is not an issue because the only required information is whether a particle passes through some region but not when it does so.
In these cases, phase errors are of no concern and the Boris algorithm is highly competitive, combining low computational cost with high quality results.
There are other applications, however, where accurate phase information is crucial.
One example are particle-wave interactions triggering Alfv\'{e}n instabilities due to resonances between orbit frequencies and wave velocities~\cite{Heidbrink2008}.
Because it is only second order accurate, the Boris method requires very small time steps in those cases, creating substantial computational cost.
In these cases, methods of order higher than two can be more efficient.

For separable Hamiltonians, the development of explicit symmetric integrators has been studied for decades~\cite{Yoshida1990}.
However, the Lorentz equations~\eqref{eq:motion_equations} give rise to a non-separable Hamiltonian, making development of higher order methods challenging, see the overview by He et al.~\cite{HeEtAl2015}.
Quandt et al. suggest a high order integrator based on a Taylor series expansion and demonstrate high convergence order for relativistic and non-relativistic test cases~\cite{QuandtEtAl2007}. 
The method needs derivatives of the electric and magnetic field, though, which may be difficult to obtain.
A recently introduced new class of methods are so-called explicit symplectic shadowed Runge-Kutta methods or ESSRK for short~\cite{Tao2016}.
They are symplectic and therefore have bounded long-term energy error.
ESSRK have been shown to be more accurate than Runge-Kutta-4 with respect to both energy and phase error but also require substantially more sub-steps.
No comparison with respect to computational efficiency seems to exist.
He at al. introduce a high-order volume preserving method based on splitting and composition of low order methods~\cite{HeEtAl2016}.
A class of symmetric multi-step methods is derived by Hairer and Lubich but not analysed with respect to computational efficiency~\cite{HairerLubich2017}.
Instead of building higher order methods, Umeda~\cite{Umeda2018} constructs a three-step version of the Boris method that can be about a factor of two faster.

\reva{Spectral deferred correction (SDC), introduced by Dutt et al. in 2000~\cite{DuttEtAl2000}, are iterative time stepping methods based on collocation.
In each time step, they perform multiple sweeps with a low order integrator (often a form of Euler method) in order to generate a higher order approximation.}
This paper presents a new high order algorithm for solving the Lorentz equations~\eqref{eq:motion_equations} called \revb{Boris-GMRES-SDC or BGSDC for short}. 
Its key advantages are that it is straightforward to implement since it heavily relies on the classical Boris scheme which will be available in almost any plasma modelling code.
Furthermore, it allows to flexibly tune the order of accuracy by simply changing runtime parameters without the need to solve equations for order conditions.
SDC also provides dense output and allows to generate a high order solution anywhere within a time step. 
We use this feature to accurately compute the turning points of particles in a magnetic mirror.
The codes used to generate the numerical examples are available for download from GitHub~\cite{code1,code2}.

\revb{BGSDC} is an extension of Boris spectral deferred corrections (Boris-SDC), introduced and tested for homogeneous electric and magnetic fields by Winkel et al.~\cite{WinkelEtAl2015}.
The present paper expands their results in multiple ways.
First, it provides a slightly simplified version of the method with almost identical performance.
Second, it integrates a GMRES-based convergence accelerator, originally introduced by Huang et al.~\cite{HuangEtAl2006} for first order case problems, with Boris-SDC.
We show that this leads to a substantial improvement in the long-term energy error.
Third, it studies the performance of \revb{BGSDC} for inhomogeneous magnetic fields, in contrast to Winkel et al. who only studied the homogeneous case.

While \revb{BGSDC} can be applied to problems where an electric field is present, we focus here \reva{on tracking fast particles in the core region of the plasma in a nuclear fusion reactor.}
\reva{There, the effect of the electric field generated from particle-interactions is small, although not totally negligible, and often ignored~\cite{AkersEtAl2002,KimEtAl2016}.
However, to include the effect of $\vect{E} \times \vect{B}$ drift on the numerical accuracy of BGSDC, we add a weak, external electric field.}
Also, to be able to quantitatively compare the accuracy of \revb{BGSDC} and the classical Boris algorithm for single trajectories, we make two simplifying assumptions.
First, we only consider test cases where the magnetic field is given by a mathematical formula (a magnetic mirror and a Solev'ev equilibrium), in contrast to a real reactor where the field is given \reva{by a numerical solution to the Grad-Shafranov equation fitted to experimental data}.
Second, we neglect the stochastic models used to capture the effect of interactions of fast ions with the plasma.
An implementation of \revb{BGSDC} into the LOCUST-GPU simulation software~\cite{LOCUST} and experiments for realistic use cases for the DIIID, JET and ITER experimental fusion reactors are ongoing work.

\paragraph{Verlet-based versus Leapfrog-based Boris integrator}
Boris-SDC relies on the classical velocity-Verlet scheme applied to~\eqref{eq:motion_equations}, which reads
\begin{subequations}\label{eq:verlet_scheme}
	\begin{align}
	\vect{x}_{n+1} &= \vect{x}_n + \Delta t \left( \vect{v}_n + \dfrac{\Delta t}{2}\vect{f}(\vect{x}_n,\vect{v}_n)\right), \\
	\vect{v}_{n+1} &= \vect{v}_n + \dfrac{\Delta t}{2} \left( \vect{f}(\vect{x}_n,\vect{v}_n) + \vect{f}(\vect{x}_{n+1},\vect{v}_{n+1})\right),
	\end{align}
\end{subequations}
with $\vect{x}_{n+1} \approx \vect{x}(t_{n+1})$, $\vect{v}_{n+1} \approx \vect{v}(t_{n+1})$ being numerical approximations of the analytical solution at some time step $t_{n+1}$.
The seemingly implicit dependence in $\vect{v}_{n+1}$ is resolved using the trick sketched in Algorithm~\ref{alg:boris} introduced by Boris in 1970~\cite{Boris1970}.
What is typically referred to as ``Boris algorithm" is the staggered Leapfrog method 
\begin{subequations}
\begin{align}
	\vect{v}_{n+1/2} &= \vect{v}_{n-1/2} + \Delta t \vect{f}(\vect{x}_n, \revb{\vect{v}_n}) \\
	\vect{x}_{n+1}   &= \vect{x}_{n} + \Delta t \vect{v}_{n+1/2}
\end{align}
\end{subequations}
which can be rewritten in ``kick-drift-kick'' form
\begin{subequations}
\label{eq:leapfrog_scheme}
\begin{align}
	\vect{v}_{n+1/2} &= \revb{\vect{v}_{n}} + \dfrac{\Delta t}{2} \vect{f}(\vect{x}_n, \vect{v}_n) \\
	\vect{x}_{n+1} &= \vect{x}_n + \Delta t \vect{v}_{n+1/2} \\
	\vect{v}_{n+1} &= \vect{v}_{n+1/2} + \dfrac{\Delta t}{2} \vect{f}(\vect{x}_{n+1}, \vect{v}_{n+1}) \label{eq:leapfrog_scheme_3}
\end{align}
\end{subequations}
where the Boris-trick is used in~\eqref{eq:leapfrog_scheme_3}.
While Velocity-Verlet~\eqref{eq:verlet_scheme} and Leapfrog~\eqref{eq:leapfrog_scheme} are similar they are not equivalent, see for example the analysis by Mazur~\cite{Mazur1997}.
\revb{In particular, in the absence of an electric field, the staggered version conserves kinetic energy exactly.}
Below, we refer to~\eqref{eq:verlet_scheme} plus the Boris trick as ``unstaggered Boris'' and to~\eqref{eq:leapfrog_scheme} with Boris trick as ``staggered Boris'' method.

While a variant of Boris-SDC can be derived based on the staggered Leapfrog method, it requires additional storage of solutions at half-points and, in tests not documented here, was not found to improve performance over the velocity-Verlet based Boris-SDC.
Substantial differences between Verlet and Leapfrog seem only to arise in simulations with very large time steps with nearly no significant digits left (phase errors well above $10^{-1}$), where staggered Boris shows better stability.
In such regimes, \revb{BGSDC} is not going to be competitive anyway so that we focus here on the simpler Verlet-based variant.
\begin{algorithm2e}[!h]
	\caption{Boris' trick \reva{for unstaggered Velocity-Verlet~\eqref{eq:verlet_scheme}}. \revc{See Birdsall and Langdon~\cite[Section 4--4]{Birdsall1985} for the geometric derivation}.}\label{alg:boris}
  	\SetKwComment{Comment}{\# }{}
	\SetCommentSty{textit}
	\SetKwInOut{Input}{input}
    \SetKwInOut{Output}{output}
    \Input{$\vect{x}_{m-1}$, $\vect{x}_m$, $\vect{v}_{m-1}$, $\Delta t$}
    \Output{$\vect{v}_{m}$ solving  $\vect{v}_m = \vect{v}_{m-1} + \Delta t \vect{E}_{m-1/2} + \Delta t \frac{\vect{v}_{m-1} + \vect{v}_m}{2} \times \vect{B}(\vect{x}_m)$}
    $\vect{E}_{m-1/2} = \frac{1}{2} \left( \vect{E}(\vect{x}_{m-1}) + \vect{E}(\vect{x}_m) \right)$ \\
    $\vect{t} = \frac{\Delta t}{2} \vect{B}(\vect{x}_{m})$ \\
    $\vect{s} = 2 \vect{t} / \left( 1 + \vect{t} \cdot \vect{t} \right)$ \\
    $\vect{v}^{-} = \vect{v}_{m-1} + \frac{\Delta t}{2} \vect{E}_{m-1/2}$ \\
    $\vect{v}^{*} = \vect{v}^{-} + \vect{v}^{-} \times \vect{t}$ \\
    $\vect{v}^{+} = \vect{v}^{-} + \vect{v}^{*} \times \vect{s}$ \\
    $\vect{v}_{m} = \vect{v}^{+} + \frac{\Delta t}{2} \vect{E}_{m-1/2}$
\end{algorithm2e}

%% file: method.tex
\section{Spectral deferred corrections}
Spectral deferred corrections~\cite{DuttEtAl2000} are based on collocation.
Therefore, we first summarise the collocation formulation of the Lorentz equations~\eqref{eq:motion_equations} before deriving the GMRES-SDC algorithm.

\subsection{Collocation}
Consider a single time step $[t_n, t_{n+1}]$.
Integrating~\eqref{eq:motion_equations} from $t_n$ to some $t_n \leq t \leq t_{n+1}$ turns them into the integral equations
\begin{subequations}
\label{eq:integral_equations}
\begin{align}
	\vect{x}(t) &= \vect{x}_0 + \int_{t_n}^{t} \vect{v}(s)~ds \\
	\vect{v}(t) &= \vect{v}_0 + \int_{t_n}^{t} \vect{f}(\vect{x}(s), \vect{v}(s))~ds
\end{align}
\end{subequations}
denoting $\vect{x}_0 = \vect{x}(t_n)$ and $\vect{v}_0 = \vect{v}(t_n)$.
The exact solution at the end of the time step can theoretically be found by inserting $t = t_{n+1}$.
In the original paper introducing Boris-SDC for homogeneous fields~\cite{WinkelEtAl2015}, the second equation is substituted into the first, resulting in double integrals over $\vect{f}$ in the equation for the position $\vect{x}$.
For the test cases studied in this paper we could not see any meaningful improvement in performance and, since the substitution leads to more complicated notation, we omit it and work directly with equations~\eqref{eq:integral_equations}.

To discretise the integral equations~\eqref{eq:integral_equations} we introduce a set of quadrature nodes $t_n \leq \tau_1 < \ldots < \tau_M \leq t_{n+1}$, set $t = t_{n+1}$ and approximate
\begin{equation}
	\int_{t_n}^{t_{n+1}} \vect{v}(s)~ds \approx \sum_{m=1}^{M} q_{m} \vect{v}_m
\end{equation}
and
\begin{equation}
	\int_{t_n}^{t_{n+1}} \vect{f}(\vect{x}(s), \vect{v}(s))~ds \approx \sum_{m=1}^{M} \revb{q_m} \vect{f}(\vect{x}_m, \vect{v}_m)
\end{equation}
with $\vect{x}_j$, $\vect{v}_j$ being approximations of $\vect{x}(\tau_j)$, $\vect{v}(\tau_j)$, that is of the analytical solution at the quadrature nodes.
Then, approximations at $t_{n+1}$ can be found from
\begin{subequations}
\label{eq:update_step}
\begin{align}
	\vect{x}_{\rm new} = \vect{x}_0 + \sum_{m=1}^{M} q_m \vect{v}_m \\
	\vect{v}_{\rm new} = \vect{v}_0 + \sum_{m=1}^{M} q_m \vect{f}(\vect{x}_m, \vect{v}_m).
\end{align}
\end{subequations}

To turn this into a usable numerical method, we require equations for the $\vect{x}_m$, $\vect{v}_m$.
Those can be obtained from discrete counterparts of~\eqref{eq:integral_equations} when setting $t = \tau_m$, for $m = 1, \ldots, M$, resulting in
\begin{subequations}
\label{eq:collocation_stages}
\begin{align}
	\vect{x}_m &= \vect{x}_0 + \sum_{j=1}^{M} q_{m,j} \vect{v}_j \approx \vect{x}_0 + \int_{t_n}^{\tau_m} \vect{v}(s)~ds \\
	\vect{v}_m &= \vect{v}_0 + \sum_{j=1}^{M} q_{m,j} \vect{f}(\vect{x}_j, \vect{v}_j) \approx \vect{v}_0 + \int_{t_n}^{\tau_m} \vect{f}(\vect{x}(s), \vect{v}(s))~ds.
\end{align}
\end{subequations}
The quadrature weights $q_{m,j}$ are given by
\begin{equation}
	q_{m, j} = \int_{t_n}^{\tau_{m}} l_{j}(s)~ds
\end{equation}
with $l_j$ being Lagrange polynomials with respect to the $\tau_m$.

Solving~\eqref{eq:collocation_stages} directly using Newton's method gives rise to a collocation method.
Collocation methods are a special type of fully implicit Runge-Kutta methods with a full Butcher tableau.
Depending on the type of quadrature nodes, they have favourable properties like symmetry (Gauss-Lobatto or Gauss-Legendre nodes)~\cite[Theorem 8.9]{HairerEtAl2002_geometric} or symplecticity (Gauss-Legendre nodes)~\cite[Theorem 16.5]{HairerEtAl1993_nonstiff} and A- and B-stability~\cite[Theorem 12.9]{HairerEtAl1996_stiff}.
However, note that even for formally symplectic implicit methods, accumulation of round-off error from the nonlinear solver can still lead to energy drift~\cite{Hairer2008}.
\subsection{Boris-SDC}
By packing the solutions  $\vect{x}_m$, $\vect{v}_m$ at the quadrature nodes into a single vector
\begin{equation}
	\vect{U} = \left( \vect{x}_1, \ldots, \vect{x}_M, \vect{v}_1, \ldots, \vect{v}_M \right)^{\rm T},
\end{equation}
the discrete collocation problem~\eqref{eq:collocation_stages} can be written as
\begin{equation}
	\label{eq:collocation}
	\vect{U} - \vect{Q} \vect{F}(\vect{U}) = \vect{U}_0.
\end{equation}
with $\vect{U}_0 = \left( \vect{x}_0, \ldots, \vect{x}_0, \vect{v}_0, \ldots, \vect{v}_0 \right)$ and
\begin{equation}
	\vect{Q} = \begin{pmatrix} \tilde{\vect{Q}} \otimes \vect{I}_{3M} & \vect{0} \\ \vect{0} & \tilde{\vect{Q}} \otimes \vect{I}_{3M} \end{pmatrix} \ \text{with} \ \tilde{\vect{Q}} = \begin{pmatrix} q_{1,1} & \ldots & q_{1,M} \\ \vdots & & \vdots \\ q_{M,1} & \ldots & q_{M,M} \end{pmatrix}.
\end{equation}
see the Appendix in Winkel et al. for details~\cite{WinkelEtAl2015}.
First, consider the case where $\vect{f}$ is linear.
For the Lorentz equations, this would be the case, for example, if $\vect{B}$ is homogeneous and $\vect{E} = 0$.
In a slight abuse of notation we write $\vect{F}$ for the matrix denoting the operator
\begin{equation}
	\vect{F}(\vect{U}) = \begin{pmatrix} \vect{v}_1 \\ \vdots \\ \vect{v}_M \\ \vect{f}(\vect{x}_1, \vect{v}_1) \\ \vdots \\ \vect{f}(\vect{x}_M, \vect{v}_M) \end{pmatrix},
\end{equation}
so that the nonlinear collocation problem~\eqref{eq:collocation} reduces to the linear system
\begin{equation}
	\label{eq:matrix_coll}
	\left( \vect{I} - \vect{Q} \vect{F} \right) \vect{U} = \vect{U}_0.
\end{equation}
\revc{One sweep of Boris-SDC can be written as}
\begin{equation}
	\label{eq:sweep_matrix_form}
	\vect{U}^{k+1} = \left( \vect{I} - \vect{Q}_{\Delta} \vect{F} \right)^{-1} \vect{U}_0 + \left[ \vect{I} - \left( \vect{I} - \vect{Q}_{\Delta} \vect{F} \right)^{-1} \left( \vect{I} - \vect{Q} \vect{F} \right) \right] \vect{U}^k.
\end{equation}
with
\begin{equation}
	\vect{Q}_{\Delta} = \begin{pmatrix} \vect{Q}_{\Delta, E} \otimes \vect{I}_{3M}& \frac{1}{2} \vect{Q}_{\Delta, E}^{(2)} \otimes \vect{I}_{3M} \\ \vect{0} & \vect{Q}_{\Delta, T} \otimes \vect{I}_{3M} \end{pmatrix}
\end{equation}
where
\begin{equation}
	\vect{Q}_{\Delta, E} = \begin{pmatrix} 0 & \ldots & 0 \\ \Delta \tau_2 & 0 & \ldots \\ \Delta \tau_2 & \Delta \tau_3 & 0 \ldots \\ &  \ldots & \\ \Delta \tau_2 & \Delta \tau_3 & \ldots \Delta \tau_M & 0 \end{pmatrix}
\end{equation}
and
\begin{equation}
	\vect{Q}_{\Delta, I} = \begin{pmatrix} \Delta \tau_1 & 0 & \ldots \\ \Delta \tau_1 & \Delta \tau_2 & 0 \ldots \\ & \ldots & \\ \Delta \tau_1 & \Delta \tau_2 & \ldots & \Delta \tau_M \end{pmatrix}
\end{equation}
and $\vect{Q}_{\Delta, T} = \frac{1}{2} \left( \vect{Q}_{\Delta, E} + \vect{Q}_{\Delta, I} \right)$ and $\vect{Q}_{\Delta, E}^{(2)} := \vect{Q}_{\Delta, E} \circ \vect{Q}_{\Delta, E}$, see again Winkel et al. for details~\cite{WinkelEtAl2015}.\footnote{\revb{As pointed out by one of the reviewers, it is also possible to use an implicit midpoint rule instead of trapezoidal rule to update the velocity. In tests not documented here, this variant of Boris-SDC showed improved long-term energy errors compared to the variant using trapezoidal rule. However, it cannot directly interpreted in the form of a preconditioned iteration given by Eq.~\eqref{eq:precond_matrix_coll}, so that it is not clear how to apply GMRES acceleration to this variant. Nevertheless, it would certainly warrant further study, in particular for problems with strong nonlinearities.}}
Iteration~\eqref{eq:sweep_matrix_form} can be understood as a Picard iteration applied to the preconditioned system
\begin{equation}
	\label{eq:precond_matrix_coll}
	\left( \vect{I} - \vect{Q}_{\Delta} \vect{F} \right)^{-1} \left( \vect{I} - \vect{Q} \vect{F} \right) \vect{U} = \left( \vect{I} - \vect{Q}_{\Delta} \vect{F} \right)^{-1} \vect{U}_0.
\end{equation}
\subsection{Boris-GMRES-SDC (BGSDC)}
For linear first order differential equations, Huang et al. showed that performing $k$ iterations of the Generalized Minimum Residual (GMRES) algorithm  on~\eqref{eq:precond_matrix_coll} often gives better results than performing $k$ standard SDC iterations~\cite{HuangEtAl2006,HuangEtAl2007}.
Here, we adopt their strategy to the second order Lorentz equations for cases where the magnetic field varies only weakly over a single time step.
Note that while we rely on a self-written GMRES implementation in the accompanying code, we verified that it gives identical results to the GMRES implementation in the SciPy library~\cite{scipy}.

GMRES does not require the matrix representing the linear system or the preconditioner to be assembled explicitly.
It only requires functions that compute $\left( \vect{I} - \vect{Q} \vect{F} \right) \vect{U}$ given some $\vect{U}$ and solve
\begin{equation}
	\label{eq:preconditioner}
	\left( \vect{I} - \vect{Q}_{\Delta} \vect{F} \right) \vect{U} = \vect{b}
\end{equation}
given some $\vect{b}$~\cite{Kelley1995}.
Applying $\vect{Q} \vect{F}$ amounts to computing the sums in~\eqref{eq:collocation_stages} for $m=1, \ldots, M$ as in the original Boris-SDC.
Systems of the form $\left( \vect{I} - \vect{Q}_{\Delta} \vect{F} \right)\vect{U} = \vect{B}$ can be solved by elimination in a sweep-like fashion.
\revc{For $M=3$ nodes, Eq.~\eqref{eq:preconditioner} becomes
\begin{align*}
	\begin{bmatrix} \vect{x}_1 \\ \vect{x}_2 \\ \vect{x}_3 \end{bmatrix} - \begin{bmatrix} 0 & 0 & 0 \\ \Delta \tau_2 \vect{I} & 0 & 0 \\ \Delta \tau_2 \vect{I} & \Delta \tau_3 \vect{I} & 0  \end{bmatrix} \begin{bmatrix} \vect{v}_{1} \\ \vect{v}_{2} \\ \vect{v}_{3} \end{bmatrix} -
	\frac{1}{2} \begin{bmatrix} 0 & 0 & 0 \\ \Delta \tau_2^2 \vect{I} & 0 & 0 \\ \Delta \tau_2^2 \vect{I} & \Delta \tau_3^2 \vect{I} & 0 \end{bmatrix} \begin{bmatrix} \vect{F}(\vect{x}_1, \vect{v}_1) \\ \vect{F}(\vect{x}_2, \vect{v}_2) \\ \vect{F}(\vect{x}_3, \vect{v}_3) \end{bmatrix} = \begin{bmatrix} \vect{b}_{1} \\ \vect{b}_{2} \\ \vect{b}_{3} \end{bmatrix} 
\end{align*}
and
\begin{align*}
	\begin{bmatrix} \vect{v}_1 \\ \vect{v}_2 \\ \vect{v}_3 \end{bmatrix} - \frac{1}{2} \begin{bmatrix}  \Delta \tau_1 \vect{I} & 0 & 0 \\ (\Delta \tau_1 + \Delta \tau_2) \vect{I} & \Delta \tau_2 \vect{I} & 0 \\ (\Delta \tau_1 + \Delta \tau_2) \vect{I} & (\Delta \tau_2 + \Delta \tau_3) \vect{I} & \Delta \tau_3 \vect{I} \end{bmatrix} \begin{bmatrix} \vect{F}(\vect{x}_1, \vect{v}_1) \\ \vect{F}(\vect{x}_2, \vect{v}_2) \\ \vect{F}(\vect{x}_3, \vect{v}_3) \end{bmatrix} = \begin{bmatrix} \vect{b}_{4} \\ \vect{b}_{5} \\ \vect{b}_{6} \end{bmatrix}.
\end{align*}
This system can be solved for $\vect{x}_{i}$, $\vect{v}_i$ by computing
\begin{align*}
	\vect{x}_1 &= \vect{b}_{1} \\
	\vect{v}_1 &= \vect{b}_{4} + \frac{1}{2} \Delta \tau_1 \vect{F}(\vect{x}_1, \vect{v}_1) \\
	\vect{x}_2 &= \vect{b}_{2} + \Delta \tau_2 \vect{v}_1 + \frac{1}{2} \Delta \tau_2^2 \vect{F}(\vect{x}_1, \vect{v}_1)\\
	\vect{v}_2 &= \vect{b}_{5} + \frac{1}{2} \left( \Delta \tau_1 + \Delta \tau_2 \right) \vect{F}(\vect{x}_1, \vect{v}_1) + \frac{1}{2} \Delta \tau_2 \vect{F}(\vect{x}_2, \vect{v}_2) \\
	\vect{x}_3 &= \vect{b}_{3} + \Delta \tau_2 \vect{v}_1 + \Delta \tau_3 \vect{v}_2 + \frac{1}{2} \Delta \tau_2^2 \vect{F}(\vect{x}_1, \vect{v}_1) + \frac{1}{2} \Delta \tau_3^2 \vect{F}(\vect{x}_2, \vect{v}_2) \\
	\vect{v}_3 &= \vect{b}_{6} + \frac{1}{2} \left( \Delta \tau_1 + \Delta \tau_2 \right) \vect{F}(\vect{x}_1, \vect{v}_1) + \frac{1}{2} \left( \Delta \tau_2 + \Delta \tau_3 \right) \vect{F}(\vect{x}_2, \vect{v}_2) + \frac{1}{2} \Delta \tau_3 \vect{F}(\vect{x}_3, \vect{v}_3)
\end{align*}
using Boris' trick to compute the velocities.
The generalisation to other values of $M$ is straightforward.
}

\subsubsection*{\revb{BGSDC} for inhomogeneous magnetic fields}
GMRES is a solver for linear systems and will not work if $\vect{f}$ and thus $\vect{F}$ are nonlinear.
In their original work, Huang et al. suggest to adopt GMRES-SDC to nonlinear problems by employing an outer Newton iteration and using GMRES-SDC as inner iteration to solve the arising linear problems.
\revc{In tests not documented in this paper we found that this approach requires too many sweeps and was not competitive for the problems studied here.}
Instead, we propose a different strategy for scenarios where $\vect{B}$ is changing slowly over the course of a time step and the nonlinearity is therefore weak.

It starts with a single sweep with standard non-staggered Boris to generate approximate values $\vect{x}^0_m$, $\vect{v}^0_m$ at all nodes.
In the notation above this is equivalent to solving
\begin{equation}
	\vect{U}^0 - \vect{Q}_{\Delta} \vect{F}(\vect{U}^0) = \vect{U}_0
\end{equation}
by block-wise elimination.
Then, we linearize the function $\vect{F}$ by setting
\begin{equation}
	\vect{F}_{\rm lin}(\vect{X}_0)(\vect{U}) = \begin{pmatrix} \vect{v}_1 \\ \vdots \\ \vect{v}_M \\  \vect{f}(\vect{x}^0_1, \vect{v}_1) \\ \vdots \\ \vect{f}(\vect{x}^0_M, \vect{v}_1) \end{pmatrix}.
\end{equation}
That is, the magnetic field applied to the velocity $\vect{v}_m$ is not $\vect{B}(\vect{x}_m)$ but $\vect{B}(\vect{x}_m^0)$ and remains fixed during the GMRES iteration.
We then apply a small number of GMRES iterations to the preconditioned linearised collocation equation
\begin{equation}
	\label{eq:linearised_collocation}
	\left( \vect{I} - \vect{Q}_{\Delta} \vect{F}_{\rm lin}(\vect{X}_0) \right)^{-1} \left( \vect{I} - \vect{Q} \vect{F}_{\rm lin}(\vect{X}_0) \right) \vect{U} = \left( \vect{I} - \vect{Q}_{\Delta} \vect{F}_{\rm lin}(\vect{X}_0) \right)^{-1} \vect{U}_0.
\end{equation}
\begin{algorithm2e}[!t]
	\caption{Single time step of \revb{BGSDC}$(k_{\textrm{gmres}},  k_{\textrm{picard}})$ for weakly nonlinear problems}
  	\SetKwComment{Comment}{\# }{}
	\SetCommentSty{textit}
	\SetKwInOut{Input}{input}
    \SetKwInOut{Output}{output}
    \Input{$\vect{x}_0$, $\vect{v}_0$}
    \Output{$\vect{x}_{\rm new}$, $\vect{v}_{\rm new}$}	
    Set $\vect{U}_0 = \left( \vect{x}_0, \ldots, \vect{x}_0, \vect{v}_0, \ldots, \vect{v}_0 \right)$ \\
    Perform a single nonlinear Boris-SDC sweep to solve $\vect{U}^0 - \vect{Q}_{\Delta} \vect{F}(\vect{U}^0) = \vect{U}_0$ \\
    $\vect{U}_{\rm lin} \leftarrow$ Perform $k_{\rm gmres}$ iterations of GMRES-SDC on the linearised collocation equation~\eqref{eq:linearised_collocation}  using $\vect{U}^0$ as starting value\\
    Perform $\revb{k_{\rm picard}}$ Picard iterations $\vect{U}^k = \vect{U}_0 + \vect{Q} \vect{F}(\vect{U}^{k-1})$ with $\vect{U}^0 = \vect{U}_{\rm lin}$ \label{line:picard} \\
    Perform update step~\eqref{eq:update_step} to compute $\vect{x}_{\rm new}$, $\vect{v}_{\rm new}$
\end{algorithm2e}

For a slowly varying magnetic field this will provide an approximation $\vect{U}_{\text{lin}}$ that is close to the solution of the nonlinear collocation problem~\eqref{eq:collocation}.
We then apply a small number of Picard iterations as sketched in Algorithm~\ref{alg:picard} using $\vect{U}_{\text{lin}}$ as starting value.
Picard iterations only require application of $\vect{Q} \vect{F}$ and do not need Boris' trick, \revc{so they are computationally cheap}.
However, they only converge for starting values that are close to the collocation solution or for small time steps.
Therefore, Picard iterations alone were not found to be competitive with either \revb{standard Boris-SDC or BGSDC}.
\revc{But for weakly nonlinear problems, the solution to the linearised collocation problem~\eqref{eq:linearised_collocation} is close to the nonlinear collocation solution~\eqref{eq:collocation} so that the output from the linearised GMRES procedure is a very accurate starting value.
Using full Boris-SDC sweeps instead of Picard iterations is also possible and, in tests not documented here, resulted in smaller errors in some cases.
We found the reduction in error is likely not significant enough to justify the higher complexity of full sweeps but leave a detailed comparison for future work.}

\begin{algorithm2e}[!t]
	\caption{Picard iteration}\label{alg:picard}
  	\SetKwComment{Comment}{\# }{}
	\SetCommentSty{textit}
	\SetKwInOut{Input}{input}
    \SetKwInOut{Output}{output}
    \Input{$\vect{U}^k = \left( \vect{x}_1^{k}, \vect{v}_1^{k}, \ldots, \vect{x}_M^{k}, \vect{v}_{M}^{k} \right)$, $\vect{x}_0$, $\vect{v}_0$}
    \Output{$\vect{U}^{k+1} = \vect{U}_0 + \vect{Q} \vect{F}(\vect{U})$}
    \For{$m = 1, \ldots, M$}{
    	$\vect{x}^{k+1}_m = \vect{x}_0 + \sum_{j=1}^{M} q_{m, j} \vect{v}_j^k$ \\
    	$\vect{v}^{k+1}_m = \vect{v}_0 + \sum_{j=1}^{M} q_{m,j} \vect{f}(\vect{x}_j^k, \vect{v}_j^k)$
    }
    $\vect{U}^{k+1} = \left( \vect{x}_1^{k+1}, \vect{v}_1^{k+1}, \ldots, \vect{x}_M^{k+1}, \vect{v}_{M}^{k+1} \right)$
\end{algorithm2e}

It was recently observed that the entries in the $\vect{Q}_{\Delta}$ matrix can be changed without \revb{losing} the sweep-like structure of SDC.
For first order problems, this allows to build more efficient sweeps resembling DIRK schemes~\cite{Weiser2014}.
In particular, one can use optimization routines to find entries for $\vect{Q}_{\Delta}$ that provide rapid convergence.
We tried to adopt this approach to second order problems but were unable to find a robust strategy that delivered improved results for a reasonably wide range of parameters.

\paragraph{Computational effort}
We use the number of evaluations of $\vect{f}$ required by each method as a proxy for computational effort $W$.
While Boris' trick requires some additional computation, in realistic simulations with experimentally given magnetic fields, evaluation of $\vect{B}(\vect{x}_m)$ dominates the computational cost because of the required interpolation.
Therefore, we count each application of Boris trick as one evaluation of $\vect{f}$, ignoring the cost of computing vector products.
Non-staggered Boris~\eqref{eq:verlet_scheme} requires one evaluation of $\vect{f}$ per time step.
Thus, its total cost when computing $N_{\rm steps}$ many time steps is simply
\begin{equation}
	W_{\rm boris} = N_{\rm steps}.
\end{equation}
In contrast, the initial predictor step in \revb{BGSDC} requires $M-1$ Boris steps and the computation of $\vect{f}(\vect{x}_0, \vect{v}_0)$ for a total of $M$ evaluations.
Computing $F_{\textrm{lin}}(\vect{X}_0)$ requires $M-1$ evaluations of $\vect{f}$ for Gauss-Lobatto nodes\footnote{We experimented with Gauss-Legendre nodes but found the resulting \revb{BGSDC} method not competitive for the studied examples. Therefore, we use Gauss-Lobatto nodes throughout the paper.}.
Because we keep the magnetic field fixed in the GMRES iterations, there is no additional cost in terms of evaluations of $\vect{f}$.
Finally, Picard iterations each require $M-1$ evaluations of $\vect{f}$ and the update step requires another $M-1$.
Therefore, the total estimated cost of \revb{BGSDC} is
\begin{equation}
	\label{eq:gmres_effort}
	W_{\rm gmres} = N_{\rm steps} \left( \underbrace{M}_{\text{predictor}} + \underbrace{(M-1)}_{\text{Compute}~\vect{F}(\vect{X}_0)} + \underbrace{(M-1) K_{\text{picard}}}_{\text{Picard iteration}} + \underbrace{M-1}_{\text{Update step~\eqref{eq:update_step}}} \right).
\end{equation}
In principle, all of those steps except the predictor can be parallelised by using $M-1$ threads to do the $\vect{f}$ evaluations for all quadrature nodes in parallel.
That would allow them to be computed in the wall clock time required for a single evaluation.
Parallelisation would reduce the cost of the algorithm to
\begin{equation}
	\label{eq:gmres_parallel_effort}
W_{\textrm{gmres}} = N_{\rm steps} \left( M + 1 + K_{\textrm{picard}} + 1 + \tau_{\textrm{overhead}} \right).
\end{equation}
Here, $\tau_{\textrm{overhead}}$ accounts for any overheads, for example from threads competing for memory bandwidth. 
There are approaches available to parallelise a full SDC sweep instead of only the Picard iteration~\cite{Speck2018} but those have not yet been adopted for second order problems.
We leave those as well as the development of an effective parallel implementation and a detailed assessment of required wall clock times for future work.

%% file: mtresults.tex
%
%
\subsection{Magnetic mirror trap}
\revb{We use simple mathematical model that has similar characteristics as a magnetic mirror trap.}
The static but non-uniform magnetic field between the coils $ \mathbf{B} = (B_x, B_y, B_z)$ has components
\begin{subequations}
\label{magneticfield}
\begin{align}
B_x =  -B_0 \dfrac{xz}{z_0^2}\\
B_y =  -B_0 \dfrac{yz}{z_0^2}\\
B_z = B_0(1-\dfrac{z^2}{z_0^2}).
\end{align}
\end{subequations}
Here, $ B_0 =  \omega_B/ \alpha$ is the magnetic field at the centre of trap, $\omega_B$ is the cyclotron frequency, $ \alpha $ is the particle's charge-to-mass ratio $ \alpha = q/m $ and $ z_0 $ the distance between coil and centre.
\revb{Note that~\eqref{magneticfield} is not a valid approximation of a mirror trap's magnetic field outside of the two coils~\cite{SimpsonEtAl2001}.}
\begin{figure}[h]
	\centering	
	\includegraphics[scale=0.7]{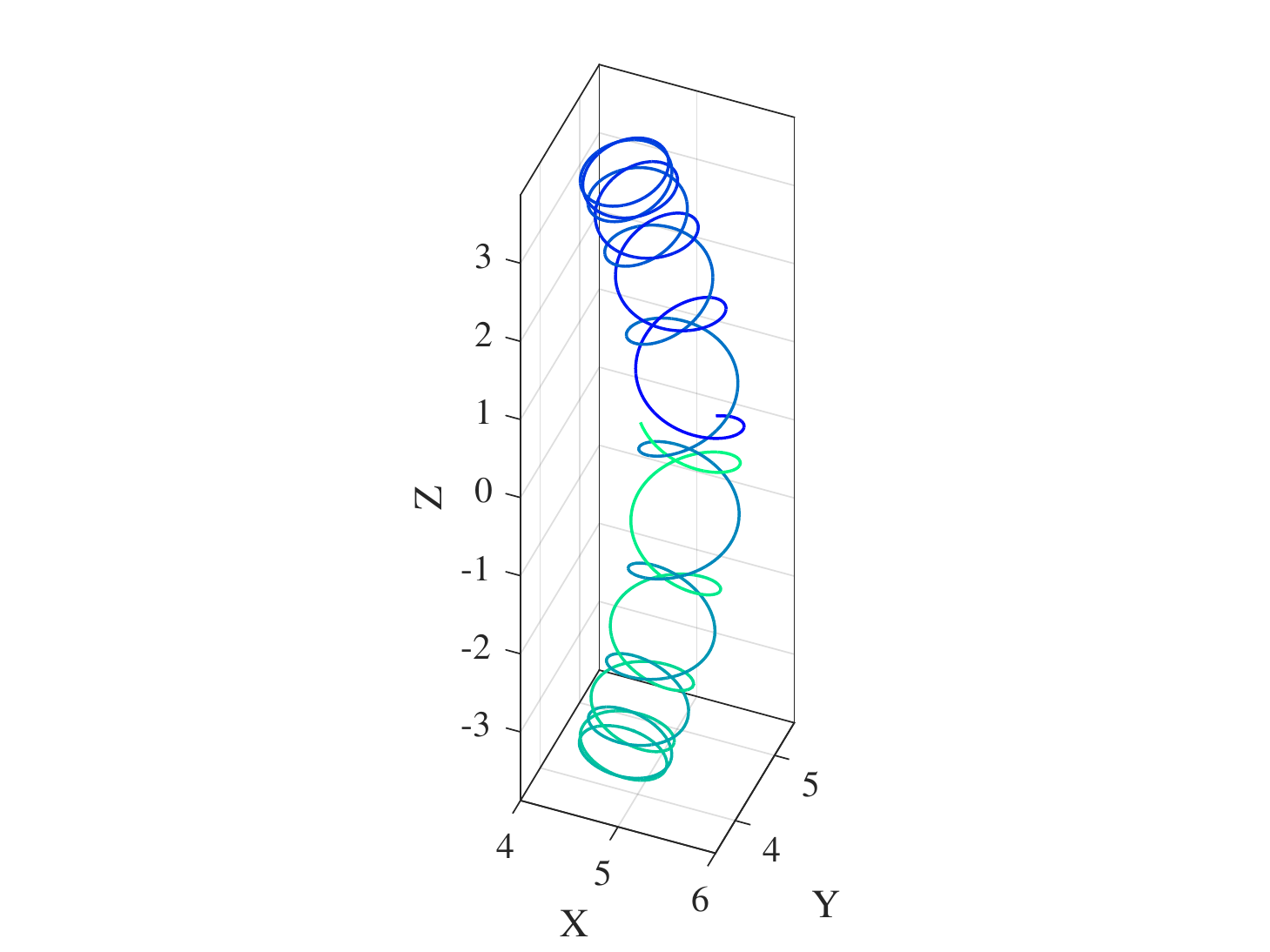}
	\caption{Example article trajectory in a magnetic mirror trap. Evolution in time is indicated by the line's changing colour from blue ($ t = 0 $) to green ($ t = t_{end} $).}
	\label{fig:trajectory}
\end{figure}

\noindent Figure~\ref{fig:trajectory} shows an example trajectory of a particle that remains vertically confined, reflecting back and forth between points at around $z = -3$ and $z = 3$.
Note that the \reva{parameters, see Table~\ref{tab:mt_parameter} (right)}, were chosen to create a recognisable trajectory and are different than the ones for \reva{the two examples reported below, see Table~\ref{tab:mt_parameter} (left and middle).}

The basic physical principle of magnetic mirroring~\cite{Budker1959, Post1958} is that charged particles in a longitudinal axially symmetric static magnetic field bounded by coils with higher value of magnetic field on both sides will be reflected from these high field side regions when moving along magnetic field lines. 
This is due to the invariance of a charged particle's magnetic moment 
\begin{equation}
	\label{eq:moment}
	\mu = \frac{1}{2}\frac{m \vect{v}_{\perp}^2}{B}, \ B = (B_x^2 + B_y^2 + B_z^2)^{1/2},
\end{equation}
in the adiabatic limit where
\begin{equation}
	\label{eq:adiabaticity}
	\varepsilon = \frac{\rho_L}{L} \ll 1 \quad \textrm{or} \quad \frac{B}{\reva{\left| \nabla B \right|}} \gg \rho_L,
\end{equation}
with $\rho_L$ being the Larmor radius of particle, $L$ the radius of curvature of the magnetic field line and the particle's velocity $\vect{v} = \vect{v}_{\perp} + \vect{v}_{||}$ being split into a part perpendicular to the magnetic field lines and a parallel part.
As a particle moves from a low field to a high field side region, $B$ increases and therefore, according to~\eqref{eq:moment}, $\vect{v}_{\perp}^2$ must increase in order to keep $\mu$ constant. 
Since the particle's kinetic energy 
\begin{equation}
	E_{\rm kin} = \frac{m \vect{v}_{\perp}^2}{2} + \dfrac{m \vect{v}_{||}^2}{2} 
\end{equation}
remains constant, the parallel velocity $\vect{v}_{||}^2$ must decrease. 
When $B$ becomes large enough, $\vect{v}_{||}$ approaches zero and the particle is reflected and travels back along the field line.

\subsubsection{\reva{Scenario 1: $\varepsilon \sim 10^{-4}$}}
%
%
\begin{table}[t]
	\caption{Parameter for \reva{scenario 1 with $\varepsilon \sim 10^{-4}$} (left), \reva{scenario 2 with $\varepsilon \sim 10^{-2}$} (middle) \reva{and for visualization (right) of a single classical particle's trajectory in Fig.~\ref{fig:trajectory} for the magnetic mirror trap.}}
	\label{tab:mt_parameter}
	\centering	
	\begin{minipage}[t]{0.32\columnwidth}
	\begin{tabular}{|ll|}
		\hline
		$ t_{end} $		& \reva{50}					\\
		$ \alpha $		& 1 					\\
		$ z_0 $			& 200						\\
		$ \omega_B $ 	& 2000  					\\
		$ \vec{x}(t=0) $& (1.0, 0.5, 0) 		\\
		$ \vec{v}(t=0) $& (100, 0, 50) 			\\
		$ N_{steps} $	& variable					\\
		$ \Delta t $	& variable			\\ \hline
	\end{tabular}
	\end{minipage}
	\begin{minipage}[t]{0.32\columnwidth}
	\begin{tabular}{|ll|}
		\hline
		$ t_{end} $		& 16					\\
		$ \alpha $		& 1 					\\
		$ z_0 $			& 16						\\
		$ \omega_B $ 	& 400  					\\
		$ \vec{x}(t=0) $& (1.0, 0, 0) 		\\
		$ \vec{v}(t=0) $& (100, 0, 50) 			\\
		$ N_{steps} $	& variable					\\
		$ \Delta t $	& variable			\\ \hline
	\end{tabular}
	\end{minipage}
	\begin{minipage}[t]{0.32\columnwidth}
	\reva{
	\begin{tabular}{|ll|}
		\hline
		$ t_{end} $		& 0.485					\\
		$ \alpha $		& 1 					\\
		$ z_0 $			& 8						\\
		$ \omega_B $ 	& 200  					\\
		$ \vec{x}(t=0) $& (5.25, 5.25, 0) 		\\
		$ \vec{v}(t=0) $& (100, 0, 50) 			\\
		$ N_{steps} $	&  1000					\\
		$ \Delta t $	& $t_{end}/N_{steps}$	\\ \hline
	\end{tabular}}
	\end{minipage}
\end{table}

In the adiabatic limit $\varepsilon \to 0$ we can determine the strength of the magnetic field at the points where the particle is reflected.
Comparing this value against the magnetic field at numerically computed reflection points allows to measure the precision of \revb{BGSDC} \reva{for very small $\varepsilon$}.
Simulation parameters are summarised in Table~\ref{tab:mt_parameter} (left) and correspond to a value $ \varepsilon \sim 8 \cdot 10^{-5}$.

Consider a particle with initial velocity $\vect{v}_0$ and position $\vect{x}_0$ and $B_0 = \left\| \vect{B}(\vect{x}_0) \right\|_2$ being the strength of the magnetic field at the particle's initial position.
Denote as $B_{\rm ref}$ the strength of the magnetic field at the reflection point and as $\vect{v}_{\perp r}$ the perpendicular velocity.
It follows from conservation of magnetic moment $\mu$ that
\begin{equation}
	\label{magmoment_conserv}
	\frac{\vect{v}_{\perp 0}^2}{B_0} = \dfrac{\vect{v}_{\perp r}^2}{B_{\rm ref}}.
\end{equation}
Conservation of kinetic energy gives
\begin{equation}
	\label{kinetic_conserv}
	\frac{m \vect{v}_{\perp 0}^2}{2} + \dfrac{m \vect{v}_{|| 0}^2}{2} = \dfrac{m\vect{v}_{\perp r}^2}{2}
\end{equation}
because, by definition,  $\vect{v}_{|| r} = 0$.
Using~\eqref{kinetic_conserv} we can substitute $\vect{v}_{\perp r} $ in~\eqref{magmoment_conserv}
\begin{equation}
	\label{eq:b_ref}
	\dfrac{B_{\rm ref}}{B_0} = \dfrac{\vect{v}_{\perp 0}^2 + \vect{v}_{|| 0}^2}{\vect{v}_{\perp 0}^2} = \dfrac{1}{\sin^2 (\varphi)},
\end{equation}
allowing us to compute $B_{\rm ref}$ directly from the initial conditions $\vect{x}_0$, $\vect{v}_0$.
Here, $ \varphi $ is the so-called pitch angle. \revb{It is assumed that the maximum magnetic field strength is on the coil $ B_{max} = \left| \mathbf{B}(z_0)\right|$ and $ \varphi_{min} = \arcsin (\sqrt{B_0/B_{max}})$.}
Note that magnetic moment is only exactly conserved in the limit $\varepsilon \to 0$.
For small but finite values of $\varepsilon$, the actual value of $B$ at the reflection point will be close to but not identical to $B_{\textrm{ref}}$.

%
%
\begin{figure}[t]
	\centering
	\includegraphics[scale=.45]{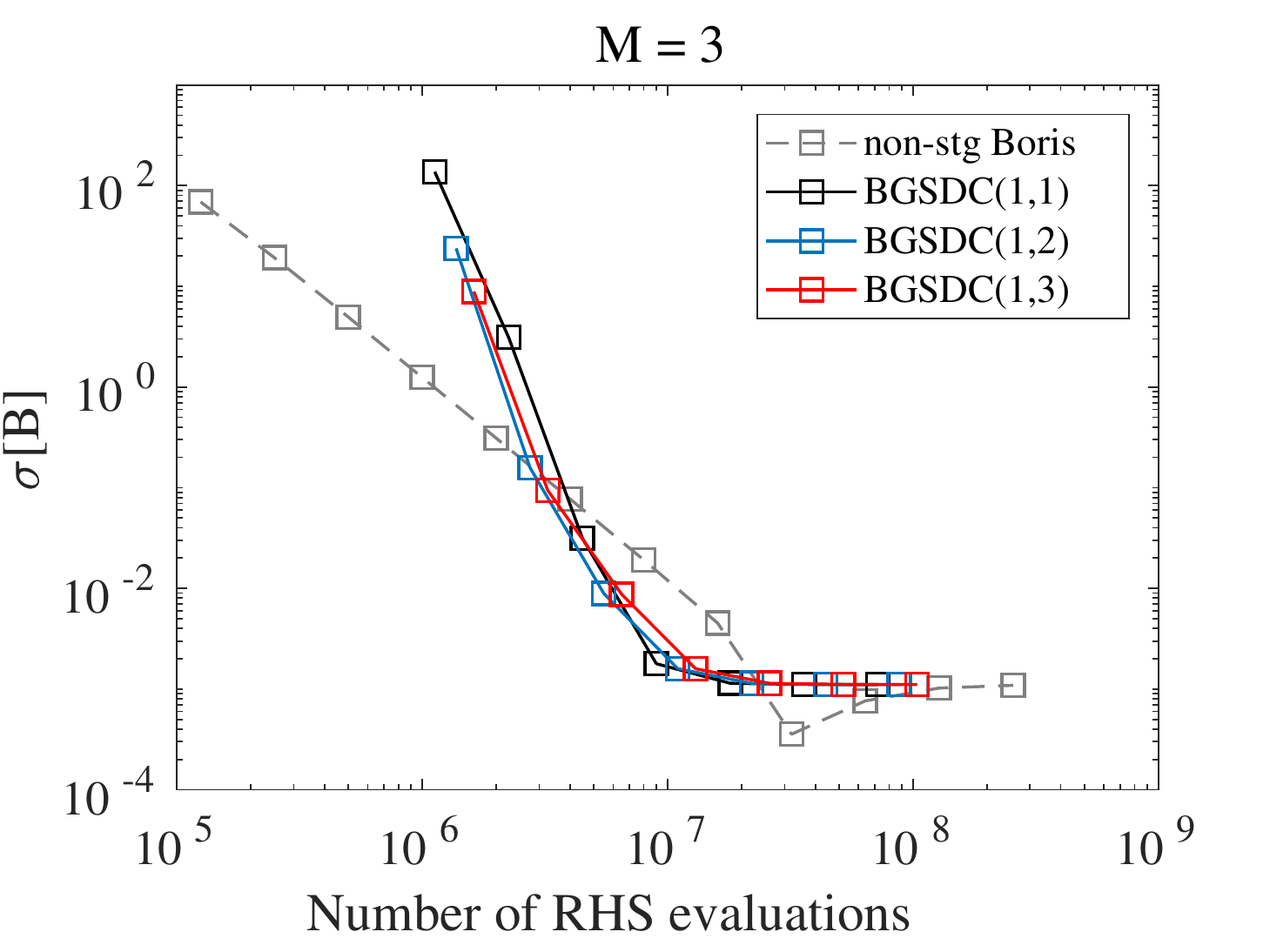}
	\includegraphics[scale=.45]{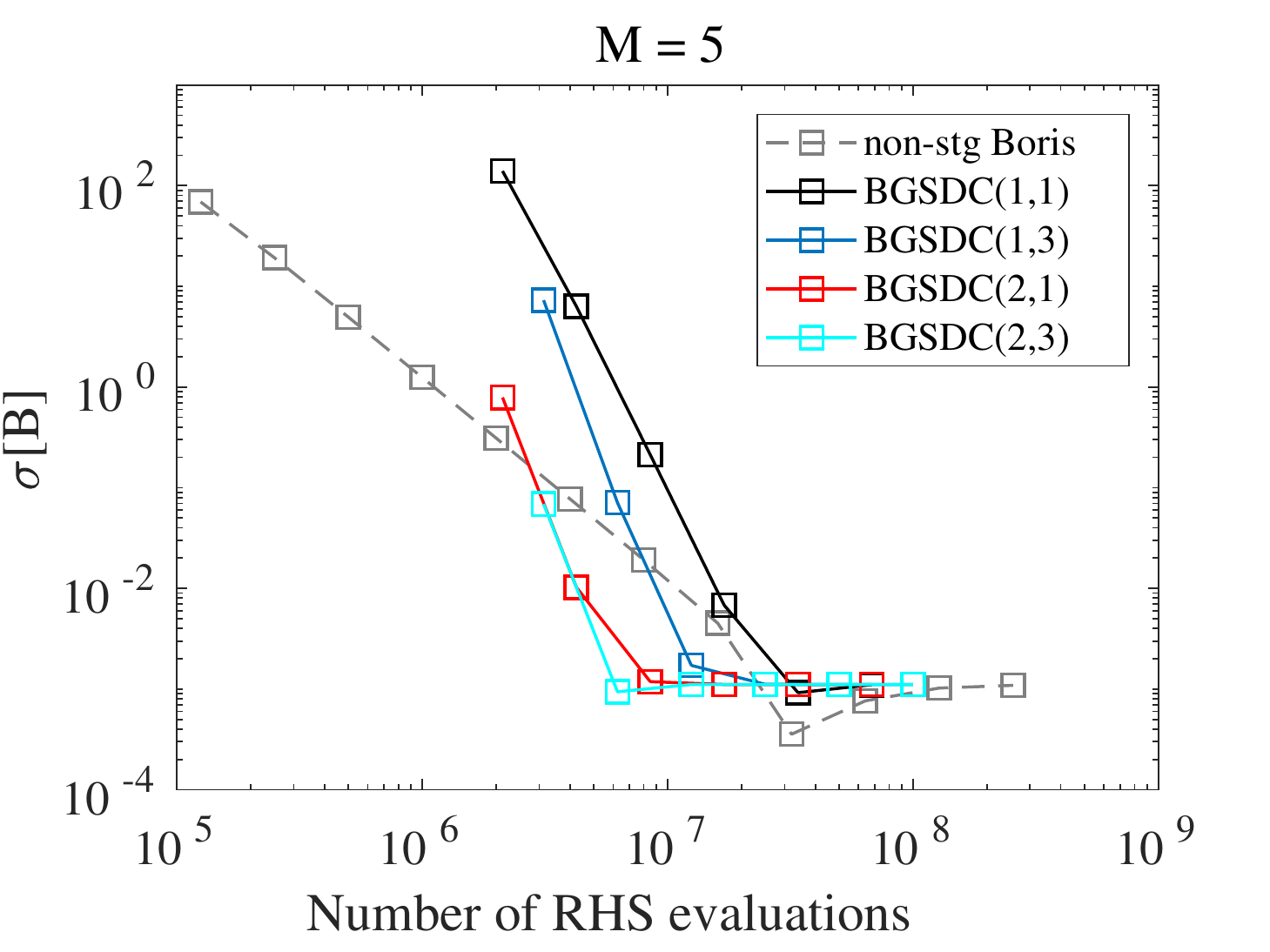}
	\caption{Error $\sigma[B]$ measured against the analytically computed value of $B$ at reflection points \reva{in the limit $\varepsilon \to 0$} plotted against the number of $\vect{f}$ evaluations.}
	\label{fig:convergence_gmth}
\end{figure}

Denote by $B_i$ the strengths of the magnetic field at the numerically computed reflection points.
We compute the $l_2$ weighted error
\begin{equation}
	\sigma [B] = \sqrt{ \dfrac{1}{N} \sum_{i=1}^{N}(B_{\rm ref} - B_i)^2},
\end{equation}
where $N$ is the number of times the particle was reflected.
To compute the $B_i$, we exploit the fact that SDC allows to reconstruct solutions at arbitrary times in a time step with high order of accuracy.
If a particle is reflected in the current time step $ [t_n, t_{n+1}]$, detected by a sign change in $v_\parallel$, we construct the Lagrange polynomial 
\begin{equation}	
	L(t) = \sum_{m=0}^{M}v_{\parallel m} \ l_m(t),
\end{equation}
using values $ \vect{x}_m, \vect{v}_m $ from intermediate nodes, where $v_{\parallel}  = v \cos(\varphi)$ and $\cos(\varphi) = \vect{v} \cdot \vect{B} / (\left\| \vect{v}\right\| \left\| \vect{B}\right\|)$.
The function $ L(t) $ interpolates $ \vect{v}_{\parallel} $ on the interval $ [t_n, t_{n+1}]$ with order $M$. 
Then, we use bisection root-finding to find the time $t_{\rm ref}$ at which $ L(t_{\rm ref}) = 0 $.
From  $t_{\rm ref}$ we can find the position $\vect{x}_{\text{ref}}$ of the reflection point using a Lagrange polynomial defined by the positions $\vect{x}_m$ at the quadrature nodes and compute the value of $B$ at that point.

Figure~\ref{fig:convergence_gmth} shows $\sigma [B]$ for \revb{BGSDC} with $M=3$ quadrature nodes (left) and with $M=5$ nodes (right) against the total number of $\vect{f}$ evaluations.
Values from the Boris method are identical in both images.
Because $B_{\rm ref}$ holds only in the adiabatic limit whereas we have a small but finite value of $\varepsilon$, errors saturate at around $10^{-3}$ for all numerically computed solutions.
For both $M=3$ and $M=5$, \revb{BGSDC}  is more efficient for precisions of $10^{-1}$ and below, requiring fewer evaluations of $\vect{f}$ than the Boris algorithm.
To reach the limit error of $10^{-3}$, \revb{BGSDC}(2,3) with $M=5$ quadrature nodes is the most efficient choice.
Boris' method requires more than ten times as many evaluations to deliver the same accuracy.

\subsubsection{\reva{Scenario 2: $\varepsilon \sim 10^{-2}$}}
\reva{Further from the adiabatic limit} we do not have an analytical solution for either the trajectory or the magnetic field at the reflection point.
Therefore, we rely on a reference solution computed numerically with a very small time step.
Simulation parameters are summarised in Table~\ref{tab:mt_parameter} (middle) and correspond to $ \varepsilon \sim 10^{-2}$.

%
%
\paragraph{Convergence order}

\begin{figure}[t]	
	\begin{subfigure}{0.5\textwidth}
		\includegraphics[scale=0.5]{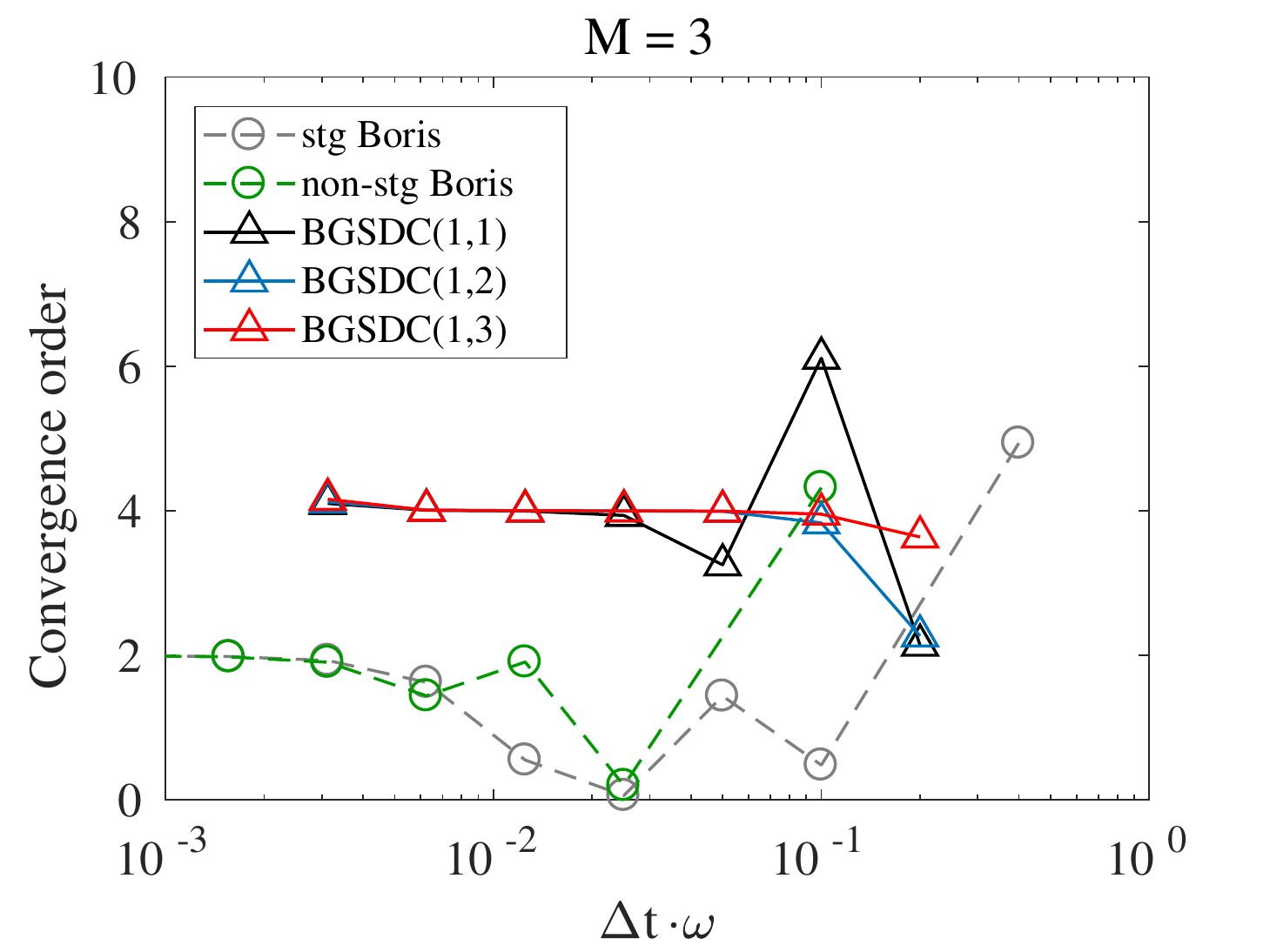} 
	\end{subfigure}
	\hfill
	\begin{subfigure}{0.5\textwidth}
		\includegraphics[scale=0.5]{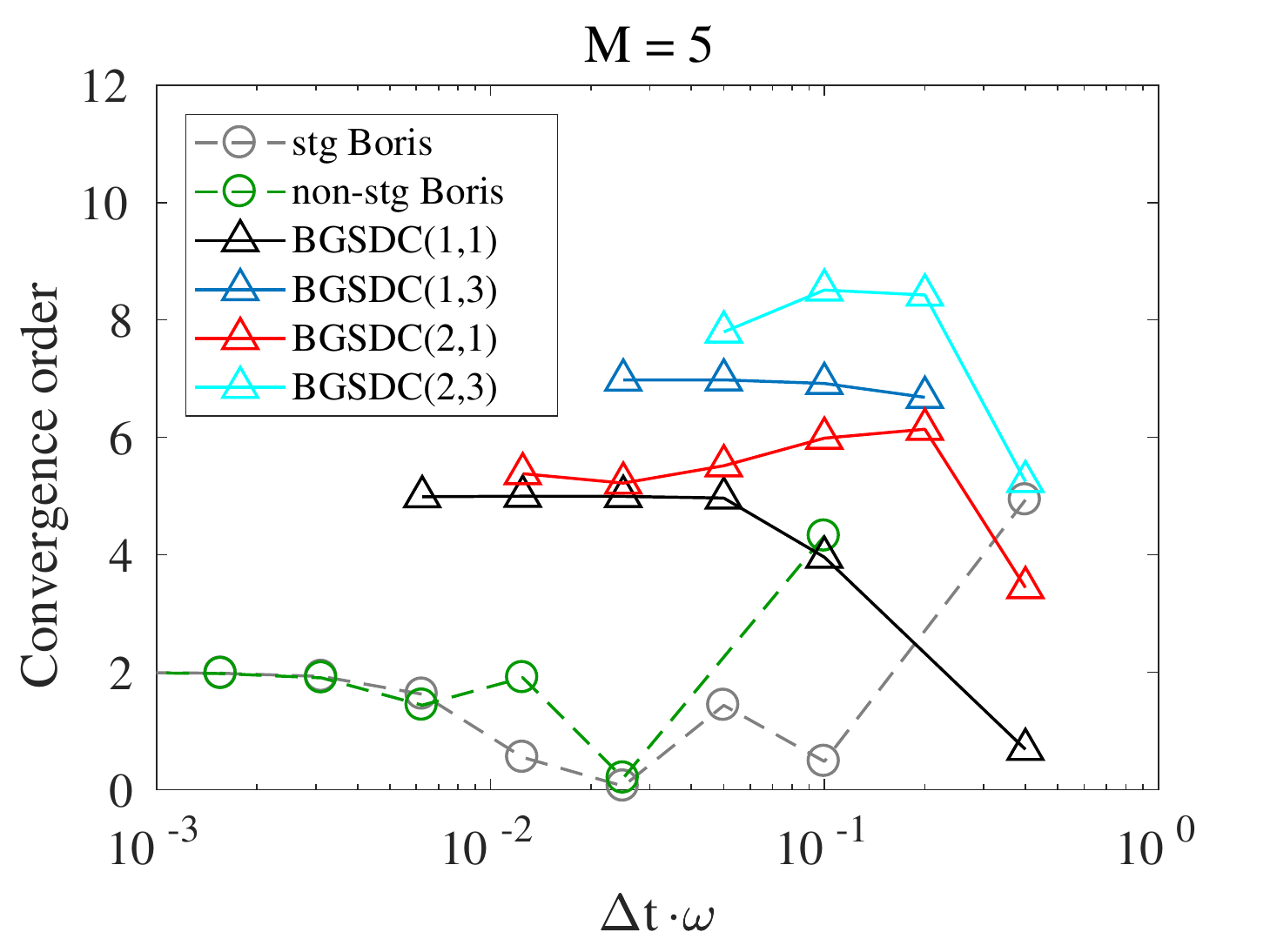}
	\end{subfigure}	
	\caption{Numerically observed convergence order $p$ for the $x$ coordinate of the particle's final position in the magnetic mirror as a function of time step size for $M=3$ (left) and $M=5$ (right)  Gauss-Lobatto collocation nodes per time step and different fixed numbers of iterations per SDC sweep.}
	\label{fig:convergencerate}
\end{figure}

\revc{Numerically computed convergence rates} for simulations with time steps from $\Delta t_0\omega = 0.0015625$ to $\Delta t_N\omega = 0.4$ are shown in Figure~\ref{fig:convergencerate} for $M=3$ nodes (left) and $M=5$ nodes (right).
While we only analyse the final error in the $x$ component of a particle's final position, results for the other position components or velocities are similar and can be generated using the published code.
Both variants of Boris achieve their theoretically expected order of $p=2$ for resolutions below $\Delta t \omega < 10^{-2}$, that is approximately $100$ steps per gyro-period.
\revb{BGSDC} with $M=3$ nodes and (1,2) and (1,3) iterations achieve the fourth order accuracy of the underlying collocation solution for $\Delta t \omega_{\rm B} < 10^{-1}$.
\revb{BGSDC}(1,1) requires a slightly smaller time step to show order $p=4$.
For $M=5$, \revb{BGSDC}(1,1) and \revb{BGSDC}(2,1) both converge with order $p \approx 5$.
This is due to having only a single Picard iteration to adjust for the nonlinearity.
Using (1,3) iterations gives order $p=7$ while (2,3) delivers the theoretical convergence order of $p=8$ of the underlying collocation solution.
Although the more complex interplay between GMRES and Picard iterations does not allow a simple heuristic like two orders per iteration that was found for non-accelerated Boris-SDC~\cite{WinkelEtAl2015}, these results show that \revb{BGSDC} can deliver high orders of convergence by changing the runtime parameter $M$ and $(K_{\textrm{gmres}}, K_{\textrm{picard}})$.

\begin{figure}[!t]	
	\begin{subfigure}{0.5\textwidth}
		\includegraphics[scale=0.50]{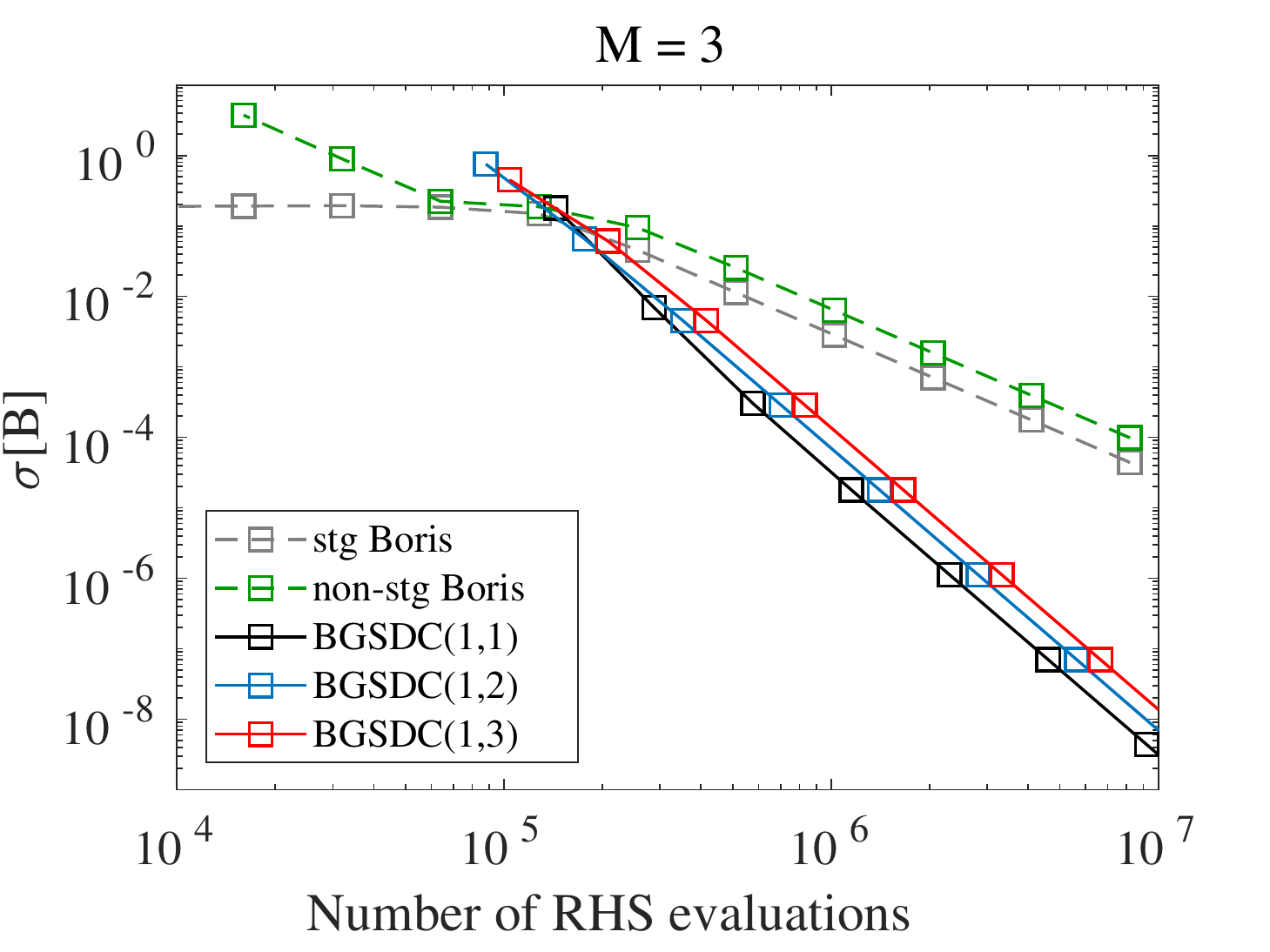} 
	\end{subfigure}
	\hfill
	\begin{subfigure}{0.5\textwidth}
		\includegraphics[scale=0.50]{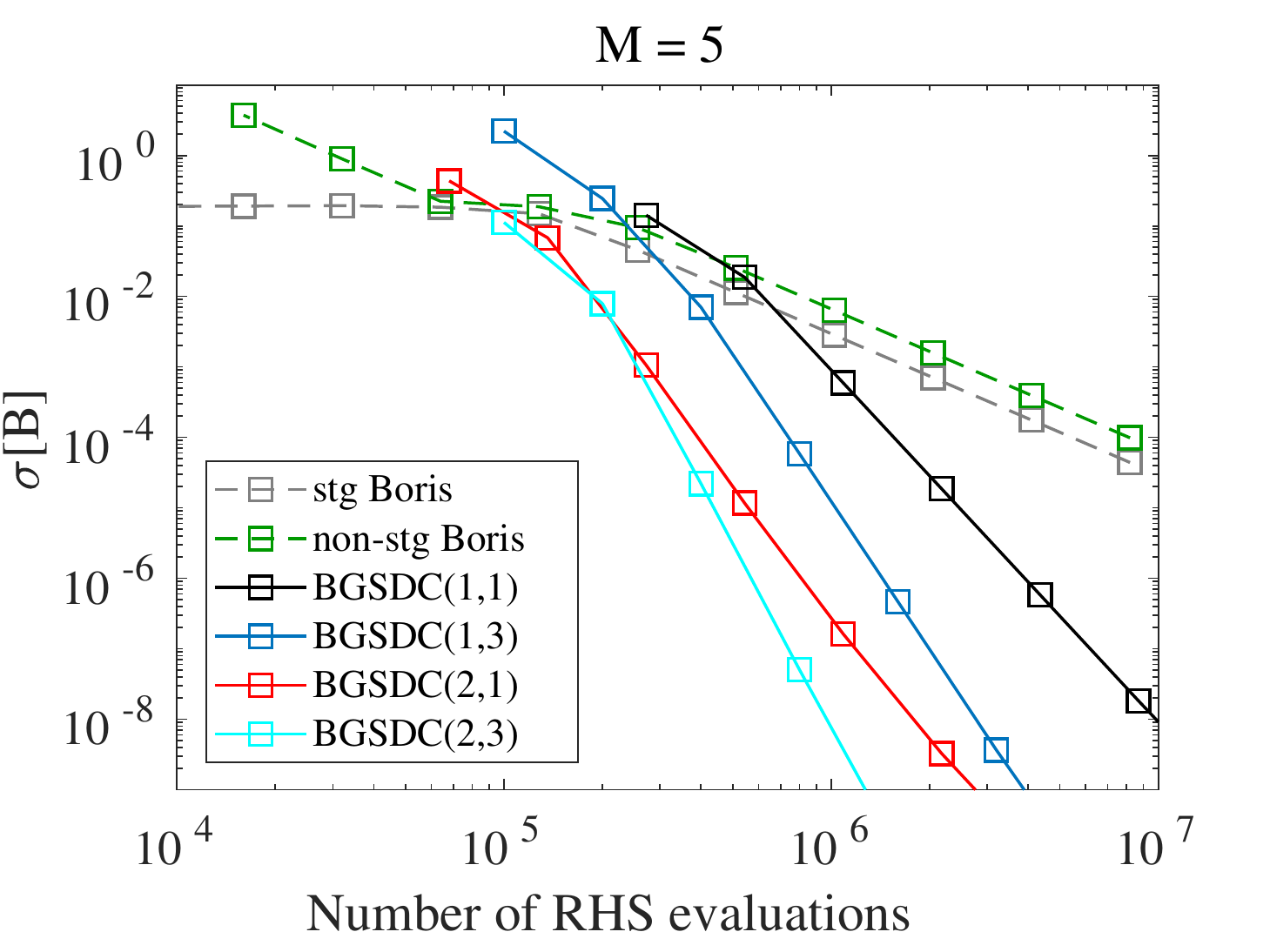}
	\end{subfigure}	
	\caption{Error $\sigma[B]$ of the magnetic field value at all reflection points during particle's motion in the magnetic mirror trap as a function of the number of r.h.s. evaluations performed for 3 and 5 Gauss-Lobatto collocation nodes per time step and different number of SDC iterations. The curves for the different runs result from varying the total number of time steps for fixed $ t_{end} $. The classical Boris integrator's (both staggered and non-staggered) convergence is shown for comparison.}
	\label{fig:evaluations}
\end{figure}

%
%
\paragraph{Work-precision analysis}
We compare the strength of the magnetic field at the reflection point against the values delivered by a reference simulation with $\Delta t\omega = 0.005$ using standard Boris-SDC with $M  = 5$ and 6 iterations.
Figure~\ref{fig:evaluations} shows the resulting error $\sigma[B]$ against the total number of $\vect{f}$ evaluations for the non-staggered Boris and \revb{BGSDC} with $M=3$ (left) and $M=5$ (right) quadrature nodes with varying numbers of iterations.
For $M=3$, all \revb{BGSDC} variants converge with order $p=4$, in line with the order of the underlying collocation method.
Increasing the number of iterations improves accuracy when keeping the time step $\Delta t$ fixed, but this does not offset the additional computational work.
Throughout, \revb{BGSDC}(1,1) is slightly more efficient than the other \revb{BGSDC} variants.
To achieve errors of $10^{-1}$ and below, \revb{BGSDC} is more efficient than Boris.
It delivers a fixed accuracy with fewer $\vect{f}$ evaluations or delivers a smaller error with the same amount of computational effort.
For errors of $10^{-3}$, the reduction in computational effort is about a factor of ten.
For $M=5$, we \revb{observe} higher convergence orders for \revb{BGSDC}(2,1) and \revb{BGSDC}(2,3), indicated by the steeper slopes.
Using (2,3) iterations delivers the most efficient method for errors below $10^{-5}$ while both \revb{BGSDC}(2,1) and \revb{BGSDC}(2,3) are about equally effective for errors up to $10^{-1}$.
Only for errors above $0.1$ does the Boris algorithm become competitive.
Note that staggered Boris gives slightly smaller errors than the non-staggered Boris, but the difference is very small.
Only for very large time steps does a substantial difference emerge.
The error for staggered Boris remains bounded at roughly $0.1$ while the error for non-staggered Boris continues to increase.
%
%

%
%
\paragraph{Long-time energy error}
Boris' method conserves phase-space volume~\cite{QinEtAl2013} which typically means a bounded long-term energy error.
For Boris-SDC and \revb{BGSDC}, depending on the choice of quadrature nodes, the collocation solution is either symmetric (Gauss-Lobatto) or symplectic (Gauss-Legendre) and will also have bounded long-term energy error, see the discussion in Winkel et al. and references therein.
However, for small numbers of iterations, both methods exhibit some energy drift.

\begin{figure}[t]	
	\begin{subfigure}{0.5\textwidth}
		\includegraphics[scale=0.5]{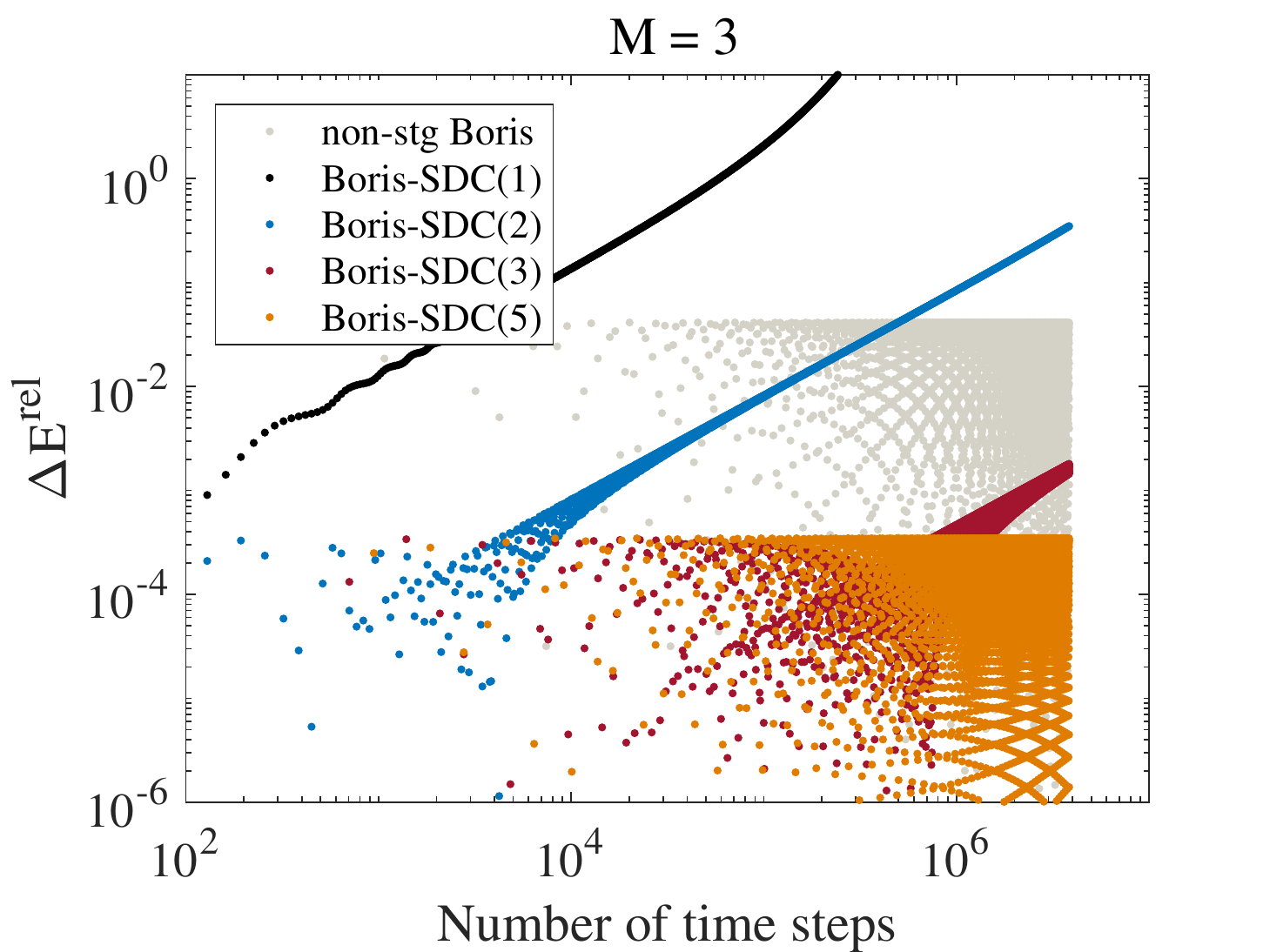} 
	\end{subfigure}
	\hfill
	\begin{subfigure}{0.5\textwidth}
		\includegraphics[scale=0.5]{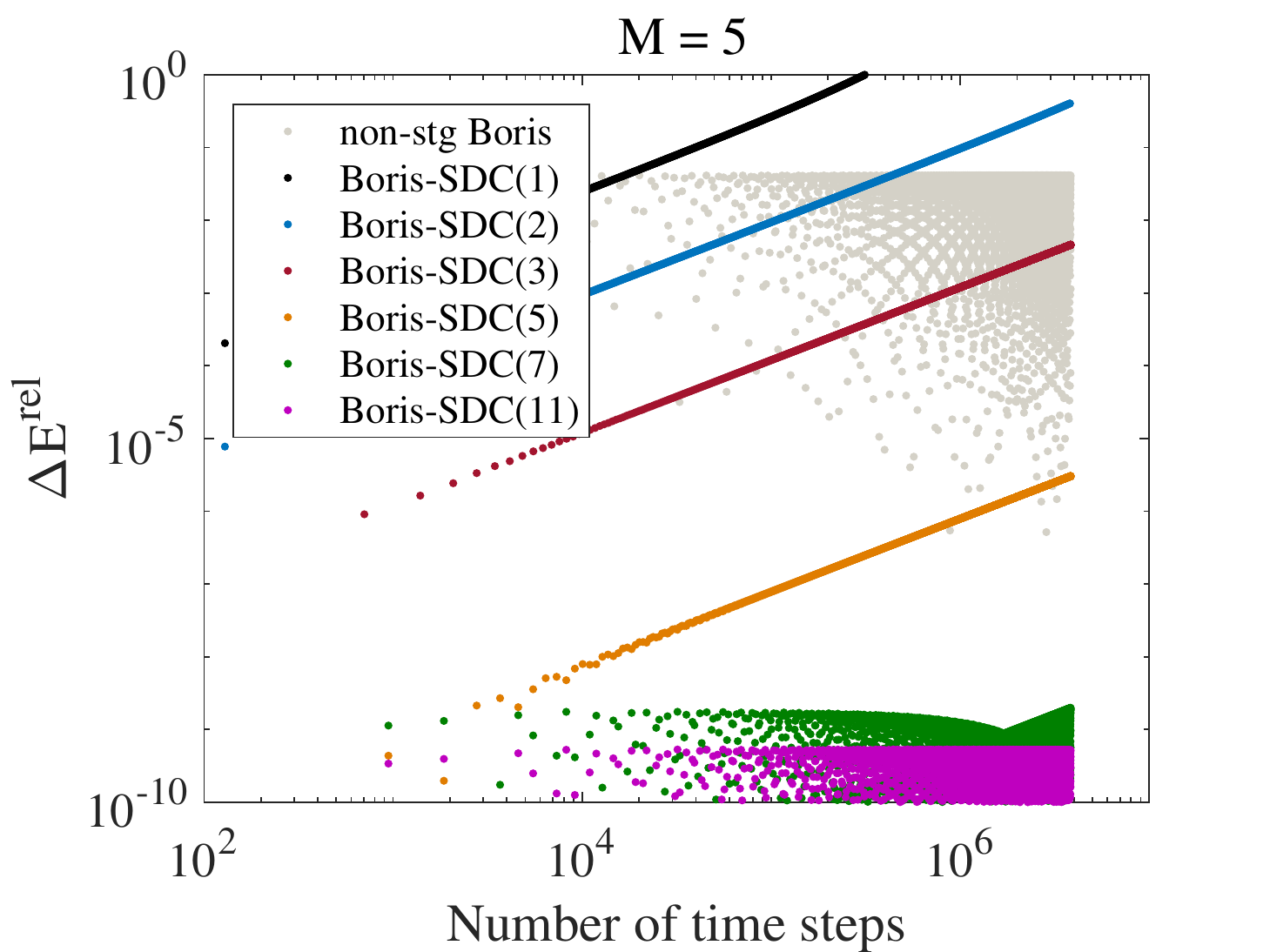}
	\end{subfigure}	
	\begin{subfigure}{0.5\textwidth}
		\includegraphics[scale=0.5]{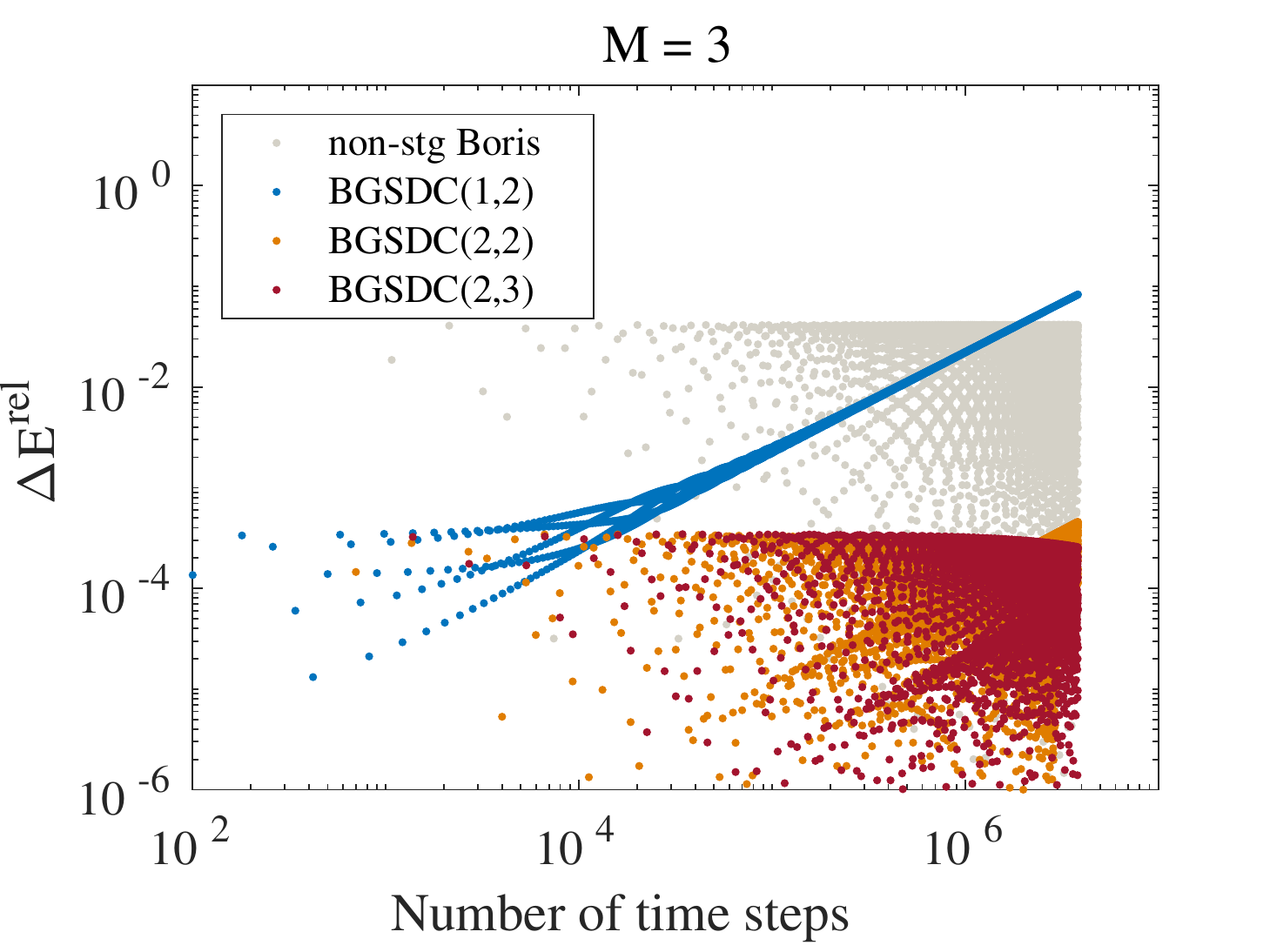} 
	\end{subfigure}
	\hfill
	\begin{subfigure}{0.5\textwidth}
		\includegraphics[scale=0.5]{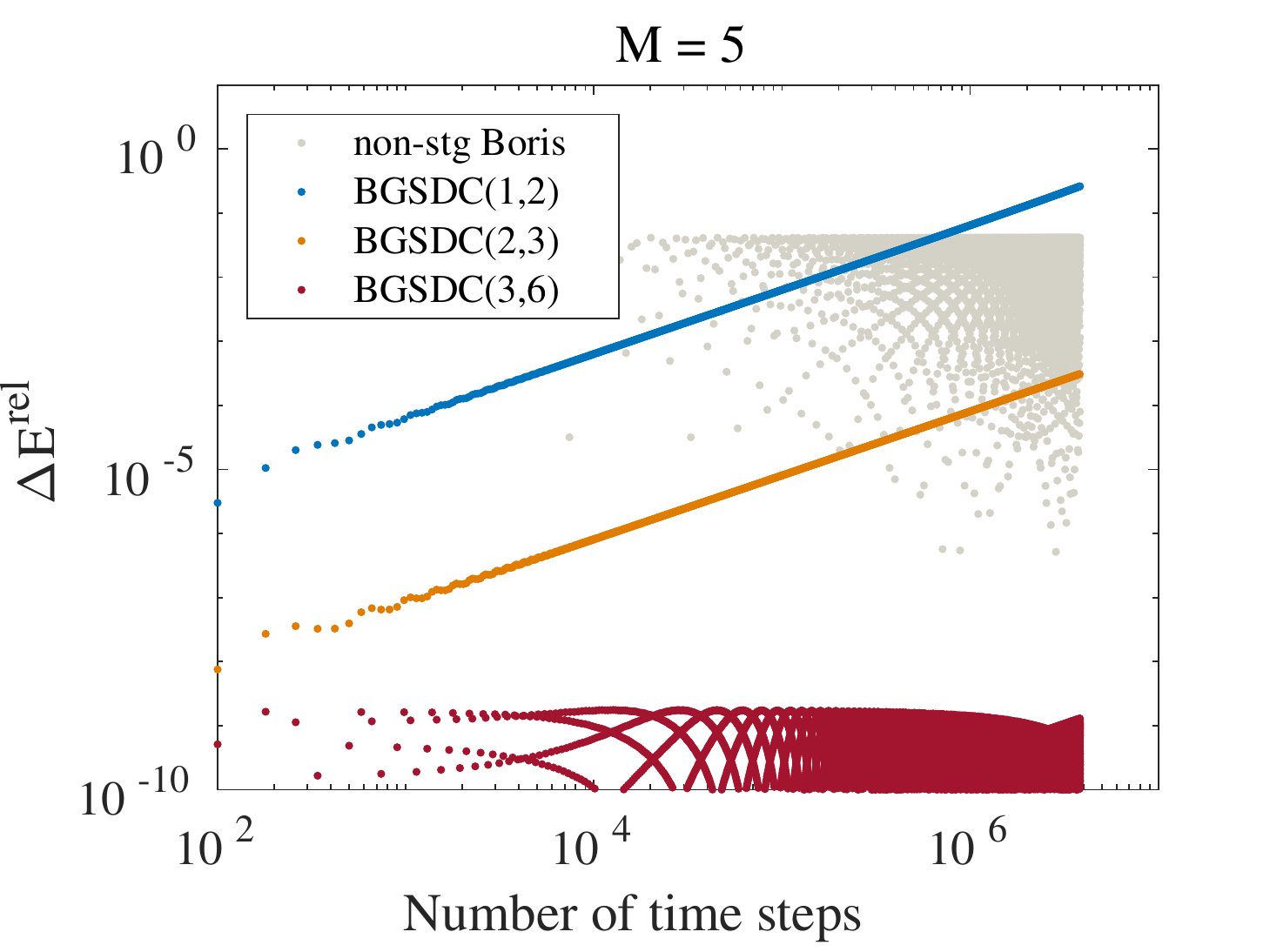}
	\end{subfigure}	
	\caption{Relative energy error \reva{for Boris-SDC (top) and \revb{BGSDC} (bottom)} over 4 million time steps with two time steps per gyro radius ($ \Delta t\omega = 0.5$) for $M=3$ (left) and $M=5$ (right) Gauss-Lobatto collocation nodes.}
	\label{fig:energy}
\end{figure}

Figure~\ref{fig:energy} shows the relative error in the total energy over $N_{steps} = 3,840,000$ time steps (with $ \Delta t\omega = 0.5$ and $t_{end} = 4800 $) for $ M = 3$ (left figures) and $M = 5$ (right figures) Gauss-Lobatto nodes and different iteration numbers.
The two upper figures show standard Boris-SDC, the lower ones \revb{BGSDC}.
Except for the larger $t_{\textrm{end}}$, parameters are identical to those used for $\varepsilon \sim 10^{-2}$, see Table~\ref{tab:mt_parameter} (middle).
As expected, \revb{non-staggered} Boris shows no drift, however its energy error is quite large at around $5 \times 10^{-2}$. 

For a small number of iterations, Boris-SDC has not yet recovered the symmetry of the underlying collocation method and shows noticeable energy drift.
However, for three iterations, after almost 4 million time steps, the energy error is still smaller than the one from Boris method for both $M = 3$ and $M = 5$ nodes. 
For five iterations and $M=3$ nodes the method has converged and the energy error remains constant.
For $M = 5$ nodes and three iterations there is drift, but the final energy error is several orders of magnitudes smaller than for Boris.
Eleven iterations are required for $M=5$ nodes for Boris-SDC to recover the bounded energy error from the collocation solution.

The lower two figures show the energy error obtained by \revb{BGSDC}.
Again, for small numbers of iteration some energy drift is observed and \revb{BGSDC}(1,2) has a larger final energy error than Boris for both $M=3$ and $M=5$.
\revb{BGSDC} requires fewer iterations than Boris-SDC to recover the bounded energy from the collocation solution.
For $M=3$, (2,3) iterations are enough (compared to five full sweeps with standard Boris-SDC) while for $M=5$ (3,6) iterations suffice, compared to 11 full sweeps.
\revc{Although Picard iterations and Boris-SDC sweeps both need $M-1$ evaluations of $\vect{f}$, Picard iterations don't require application of the preconditioner $\left( \vect{I} - \vect{Q}_{\Delta} \vect{F} \right)$ and will thus be computationally cheaper in terms of runtime.}
Therefore, \revb{BGSDC} delivers a smaller energy error for less computational work than Boris-SDC.

%% file: solresults.tex
%
%
\subsection{Solev'ev equilibrium}
As a second test case we consider the Solev'ev equilibrium~\cite{Solovev1968,Zheng1996} \revb{with an added simple radial electric field}. The magnetic field is given by
\begin{equation}
	\label{Solov'ev}
	\begin{cases}
	B_R = -(\dfrac{2 \tilde{y}}{\sigma^2})(1 - 0.25 \epsilon^2) \dfrac{(1 + \kappa \epsilon \tilde{x}(2 + \epsilon \tilde{x}))}{\psi R} \\
	B_Z = 4(1 + \epsilon \tilde{x}) \dfrac{(\tilde{x} - 0.5 \epsilon (1 - \tilde{x}^2) + (1 - 0.25 \epsilon^2) \tilde{y}^2 \kappa \dfrac{\epsilon}{\sigma^2})}{ \psi R ((r_{ma} - r_{mi})/z_0)}\\
	B_\phi = \dfrac{B_\phi^0}{R}
	\end{cases}.
\end{equation}
Here, $ \sigma, \epsilon, \kappa, \psi,  r_{ma}, r_{mi}, z_m, B_\phi^0$ are constants given in Table~\ref{tab:solevev}, chosen to model an equilibrium similar to the one in the Joint European Torus (JET) fusion reactor\footnote{We thank Dr Rob Akers from Culham Centre for Fusion Energy for providing this test case and the parameter to model the JET equilibrium.}.
Furthermore, $ \tilde{x}, \tilde{y} $ are the intermediate coordinates
\begin{equation}
\begin{gathered}
	\tilde{x} = 2 \dfrac{(R - r_{mi})}{(r_{ma} - r_{rmi})} - 1, \\
	\tilde{y} = (Z - z_m)/z_0.
\end{gathered}
\end{equation}
\revb{The radial electric field $ E_r = E_0 r^2/r_a^2$ is given in toroidal coordinates $ (r, \phi, \theta) $ which are connected to a cylindrical system via $ R = R_0 + r\cos(\theta), Z = Z_0 + r\sin(\theta) $. 
Constants $ E_0, r_a, R_0, Z_0 $ are given in Table~\ref{tab:solevev}.}

\begin{figure}[t]
	\centering	
	\includegraphics[width=0.49\textwidth, height = 5cm, keepaspectratio]{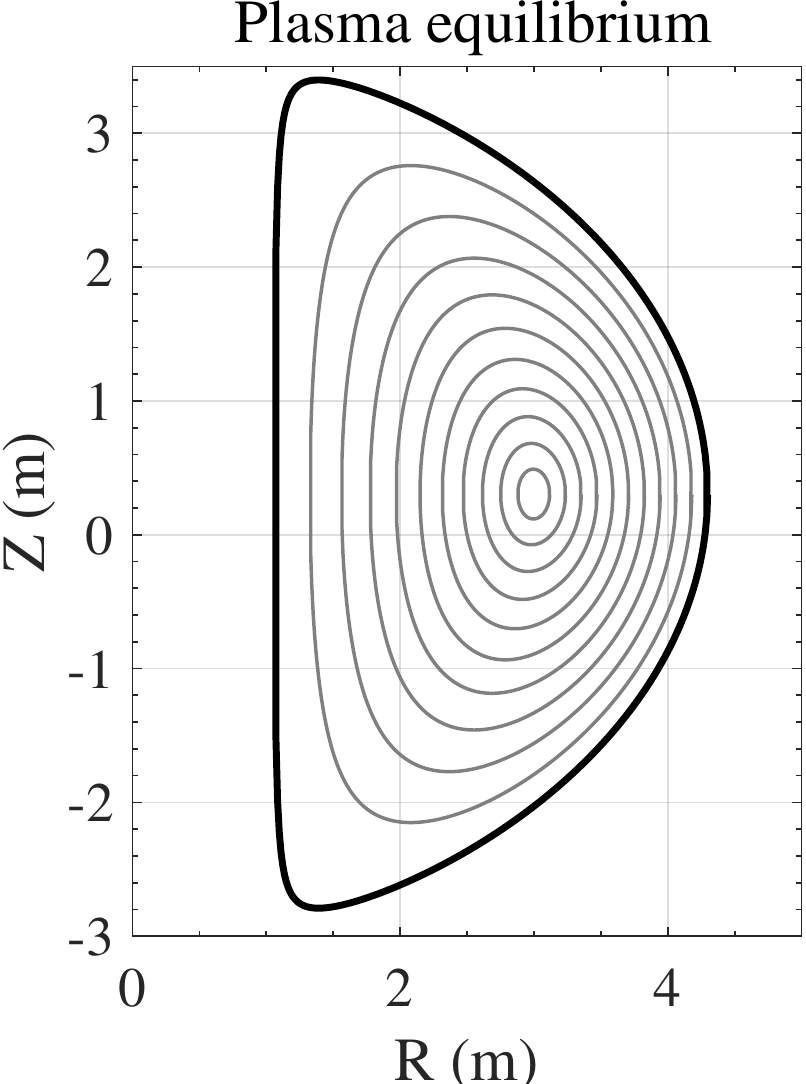}
	\includegraphics[width=0.49\textwidth, height = 5cm, keepaspectratio]{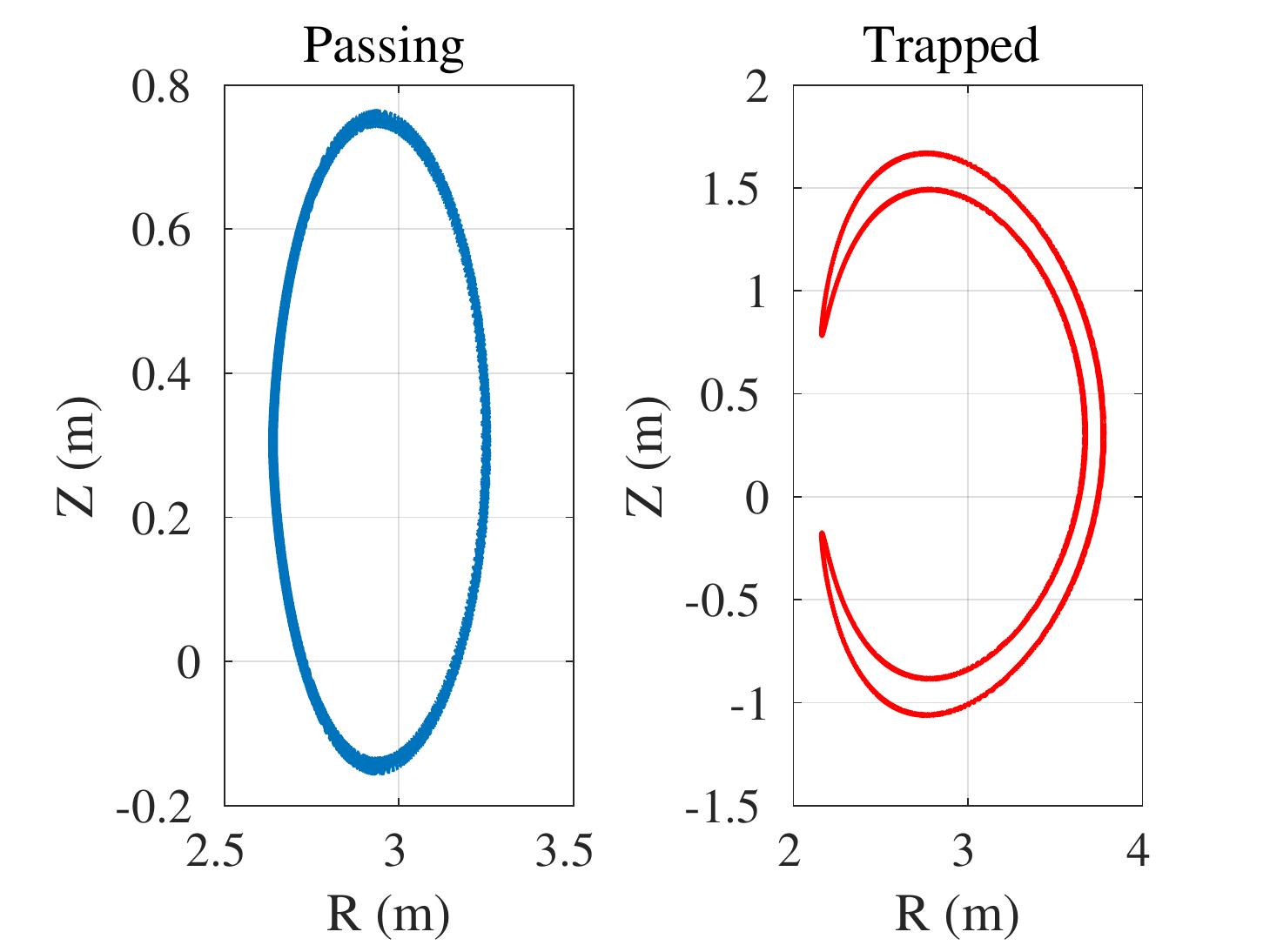}	
	\caption{Magnetic surfaces for JET-like equilibria with parameters from Table 3 (left) and example trajectories of a passing (middle) and trapped particle (right).}
	\label{fig:soleq}
\end{figure}

Fig.~\ref{fig:soleq} shows the magnetic surfaces in a $R$-$Z$ cross-section as well as the two example trajectories studied below.
One is for a passing particle that continues to perform full revolutions in the reactor's magnetic field.
The second is a trapped particle which changes direction at some point of its orbit.
It thus fails to complete a full revolution and instead travels on a so-called ``banana-orbit''.
Initial position and velocity for both the passing and the trapped particle are given in Table~\ref{tab:solevev_initial}.
\begin{table}[t]
	\centering
	\begin{tabular}{ll}
		\hline
		$ \sigma $		&  \ 1.46387369075	\\
		$ \epsilon $	&  \ 0.22615668214	\\ 
		$ \kappa $		&  \ 1.43320389205	\\ 
		$ \psi $		&  \ 1.13333149039 $\left[ \tesla^{-1}\meter^{-1} \right] $	\\ 
		$ r_{ma} $		&  \ 3.83120489000 $\left[ \meter\right] $  \\
		$ r_{mi} $		&  \ 1.96085203000 $\left[ \meter\right] $ 	\\ 
		$ z_m $			&  \ 0.30397316800 $\left[ \meter\right] $ 	\\
		$ z_0 $			&  \ 1.0 $\left[ \meter\right] $ 	\\  
		$ B_\phi^0 $	& -9.96056843000 $\left[ \tesla \meter \right] $	\\  
		\revb{$ E_0 $}			& \	50000 $\left[ \volt /  \meter \right]$ \\  		
		\revb{$ r_a $}			& \	1.5 $\left[ \meter\right] $\\ 
		\revb{$ R_0 $}			& \	3.00045800 $\left[ \meter\right] $\\ 
		\revb{$ Z_0 $}			& \	0.30397317 $\left[ \meter\right] $ \\ \hline
	\end{tabular}
	\caption{Parameters needed to reconstruct the magnetic field.}
	\label{tab:solevev}
\end{table}
\begin{table}[th]
		\centering
	\begin{tabular}{lll}
			\hline
			&  \textit{Passing} & \textit{Trapped} \\
			\hline
	$ x $	&\ 2.1889641172761  &\ 3.0852639552352  \\
	$ y $	&\ 0  				&\ 0  					\\
	$ z $	&\ 0.8635434778595  &-0.0732997600262  \\
	$ v_x $	&\ 2269604.3143406  &\ 814158.31065935   \\
	$ v_y $	&\ 292264.06108651  &\ 931354.18390575   \\
	$ v_z $	&-338526.06660893  	&\ 1793580.5493877 \\	\hline 
	\end{tabular}
	\caption{Initial position and velocity for a passing and trapped particle in the Solev'ev equilibrium. The charge-to-mass ratio is $\alpha = 47918787.60368$, \revb{gyrofrequency ${ \omega_0 = 0.159\cdot10^{-9} }$ on the magnetic axis $ (R_0,Z_0) $} and their trajectories are simulated until $t_{\rm end} = \si{10}{\ \milli\second}$.}
	\label{tab:solevev_initial}
\end{table}

%
%
\subsubsection{Accuracy}
To assess accuracy of Boris and \revb{BGSDC}, we compare their particle trajectories against a reference trajectory computed using \revb{BGSDC}(2,4) $ M=5 $ with a time step of $\Delta t = \SI{0.1}{\nano\second}$.
We choose time steps such that the time points in every run are a subset of the time points of the reference to avoid the need for interpolation.
The maximum defect in $x$ at all points of the computed trajectory reads
\begin{equation}
	\label{eq:trajectory_error}
	d_x := \max_{n = 0, \ldots, N} \left| x_n - x^{(\text{ref})}(t_n) \right|,
\end{equation}
with $t_n$, $n=0, \ldots, N$ being the time steps for the current resolution and $x^{(ref)}(t_n)$ the reference solution at those points.
Analogous expressions are used to compute $d_y$ and $d_z$ and we then take the maximum.
Note that~\eqref{eq:trajectory_error} compares positions at a specific time so that $d_{\textrm{max}}$ measures not only particle drift but also errors in phase.
Table~\ref{tab:solevev_accuracy} shows the resulting trajectory errors for the passing particle (upper two) and trapped particle (lower two) for staggered Boris and \revb{BGSDC} with $M=3$ and $M=5$ nodes.

\begin{table*}[!th]
	\centering
	\begin{minipage}[t]{\columnwidth}
	\centering
	\revb{\begin{tabular}{|c|c|c|c|} \hline
		\multicolumn{4}{|c|}{Trajectory error for passing particle for $M=3$ nodes.} \\ \hline
		$\Delta t$             & Staggered Boris & \revb{BGSDC}(1,3) & \revb{BGSDC}(2,6) \\ \hline
		\SI{0.1}{\nano\second} & 0.4362	&    -      	&  -	\\ 
		\SI{0.2}{\nano\second} & 1.7210	&    -      	&  -	\\ 
		\SI{0.5}{\nano\second} & 6.3286	& 0.2876        & 0.0661 \\ 
		\SI{1}{\nano\second}   & 6.3089	& 6.3189        & 0.9593 \\
		\SI{2}{\nano\second}   & 6.4776	& 6.4562        & 4.3799 \\ \hline    
	\end{tabular}}
	\end{minipage}\vspace*{0.5\baselineskip}
	\begin{minipage}[t]{\columnwidth}
	\centering	
	\revb{\begin{tabular}{|c|c|c|c|} \hline
		\multicolumn{4}{|c|}{Trajectory error for passing particle for $M=5$ nodes.} \\ \hline
		$\Delta t$             & Staggered Boris & \revb{BGSDC}(1,4) & \revb{BGSDC}(2,6) \\ \hline
		\SI{0.1}{\nano\second} & 0.4362	&    -      	&  -	\\ 
		\SI{0.2}{\nano\second} & 1.7210	&    -      	&  -	\\ 
		\SI{0.5}{\nano\second} & 6.3286 & 0.02351848	& 0.00000049 \\ 
		\SI{1}{\nano\second}   & 6.3089 & 2.03153722	& 0.00158861  \\
		\SI{2}{\nano\second}   & 6.4776	& 6.49396896 	& 0.42631203 \\ \hline       
	\end{tabular}}
	\end{minipage}\vspace*{0.5\baselineskip}
	\begin{minipage}[t]{\columnwidth}
	\centering	
	\revb{\begin{tabular}{|c|c|c|c|} \hline
		\multicolumn{4}{|c|}{Trajectory error for trapped particle for $M=3$ nodes.} \\ \hline
		$\Delta t$             & Staggered Boris & \revb{BGSDC}(1,3) & \revb{BGSDC}(2,6) \\ \hline
		\SI{0.1}{\nano\second} & 0.0276			&    -      	&  -	\\ 
		\SI{0.2}{\nano\second} & 0.0353			&    -      	&  -	\\ 
		\SI{0.5}{\nano\second} & 0.1721			& 0.1212	& 0.0027	\\ 
		\SI{1}{\nano\second}   & 0.6976			& 6.7962	& 0.0275 	\\
		\SI{2}{\nano\second}   & 2.7529    		& 7.3273    & 5.7878 	\\ \hline        
	\end{tabular}}
	\end{minipage}\vspace*{0.5\baselineskip}
	\begin{minipage}[t]{\columnwidth}
	\centering	
	\revb{\begin{tabular}{|c|c|c|c|} \hline
		\multicolumn{4}{|c|}{Trajectory error for trapped particle for $M=5$ nodes.} \\ \hline
		$\Delta t$             & Staggered Boris & \revb{BGSDC}(1,4) & \revb{BGSDC}(2,6) \\ \hline
		\SI{0.1}{\nano\second} & 0.0276			&    -      	&  -	\\ 
		\SI{0.2}{\nano\second} & 0.0353			&    -      	&  -	\\ 
		\SI{0.5}{\nano\second} & 0.1721			& 0.01407600   	& 0.00000022    \\ 
		\SI{1}{\nano\second}   & 0.6976			& 1.07024660 	& 0.00109417 	\\
		\SI{2}{\nano\second}   & 2.7529    		& 6.96382918	& 0.53245320 	\\ \hline     
	\end{tabular}}
	\end{minipage}	
	\caption{Maximum  deviation from reference trajectory $d_{\rm max}$ in \SI{}{\metre} for a passing (upper two) and trapped (lower two) particle in a Solev'ev equilibrium. Note that the values for the staggered Boris method in the two upper and two lower tables are identical, because the number of nodes $M$ does not affect it.}
	\label{tab:solevev_accuracy}
\end{table*}

\paragraph{Passing particle} \revb{For a passing particle, if precision in the range of millimetres is required, BGSDC(2,6) with $M=5$ nodes can deliver this with a \SI{1}{\nano\second} time step.
In contrast, staggered Boris has an error of around \SI{40}{\centi\metre} even with a \SI{0.1}{\nano\second} step.
Assuming it is converging with its theoretical order of two, staggered Boris would require a time step of $0.1 / \sqrt{0.4362/0.00159}$ or approximately \SI{0.006}{\nano\second} to be as accurate as BGSDC.
Therefore, it requires $\SI{1}{\nano\second} / \SI{0.006}{\nano\second} \approx 167$ times as many steps as BGSDC to deliver the same accuracy.
However, BGSDC(2,6) with $M=5$ nodes requires $5 + 4 + 4*6 + 4 = 37$ $\vect{f}$-evaluations per time step according to~\eqref{eq:gmres_effort} whereas Boris needs only one.
Nevertheless, BGSDC(2,6) with $\Delta t = \SI{1}{\nano\second}$ will be around $167/37 \approx 4.5$ times faster than staggered Boris with $\Delta t = \SI{0.006}{\nano\second}$ while delivering the same accuracy.}

\revb{If only centimetre precision is required, \revb{BGSDC} without parallelisation will struggle to be faster than staggered Boris.
BGSDC(2,6) with $M=3$ nodes achieves an error of about \SI{6.6}{\centi\metre} for a time step of \SI{0.5}{\nano\second}.
For the same accuracy, staggered Boris would require a time step of around $\SI{0.1}{\nano\second} / \sqrt{0.4362/0.0661} \approx \SI{0.04}{\nano\second}$.
Therefore, staggered Boris needs about $13$ times as many steps as BGSDC, but BGSDC would be around $3 + 2 + 2*6 +2 = 19$ times more expensive per step, thus making it slower.
The parallel \revb{BGSDC}(1,3) with workload model~\eqref{eq:gmres_parallel_effort} would be competitive as it is only $3 + 1 + 3 + 1 = 8$ times more expensive per step.}

\paragraph{Trapped particle} \revb{BGSDC(2,6) with $M=5$ nodes can deliver micrometre precision with a \SI{1}{\nano\second} time step.
Staggered Boris has an error of \SI{2.76}{\centi\metre} for a \SI{0.1}{\nano\second} step and would require approximately a $\SI{0.1}{\nano\second} / \sqrt{0.0276/0.0011} \approx \SI{0.02}{\nano\second}$ time step to be comparable in precision to BGSDC.
Note that we again assume that staggered Boris converges with its full second order accuracy, even though the reduction in error from \SI{0.2}{\nano\second} to \SI{0.1}{\nano\second} time step suggests that this is not yet the case.
Staggered Boris therefore needs at least $1 / 0.02 = 50$ times more steps than BGSDC while every step of BGSDC(2,6) is 37 times more expensive.
Therefore, we expect BGSDC to be at least $50/37 \approx 1.3$ times faster than staggered Boris.}

\revb{Centimetre precision can be delivered by BGSDC(2,6) with $M=3$ nodes and a \SI{1}{\nano\second} time step or by staggered Boris with a time step of \SI{0.1}{\nano\second}.
Thus, staggered Boris requires only about 10 times as many steps but BGSDC(2,6) is 19 times more expensive per step, making it slower.
Even the parallel BGSDC versions is 11 times more expensive so that additional improvements are required for it to be competitive.}

%
%
\subsubsection{Long-time energy error}
\begin{figure}[t]
	\begin{subfigure}{0.5\textwidth}
		\includegraphics[scale=0.5]{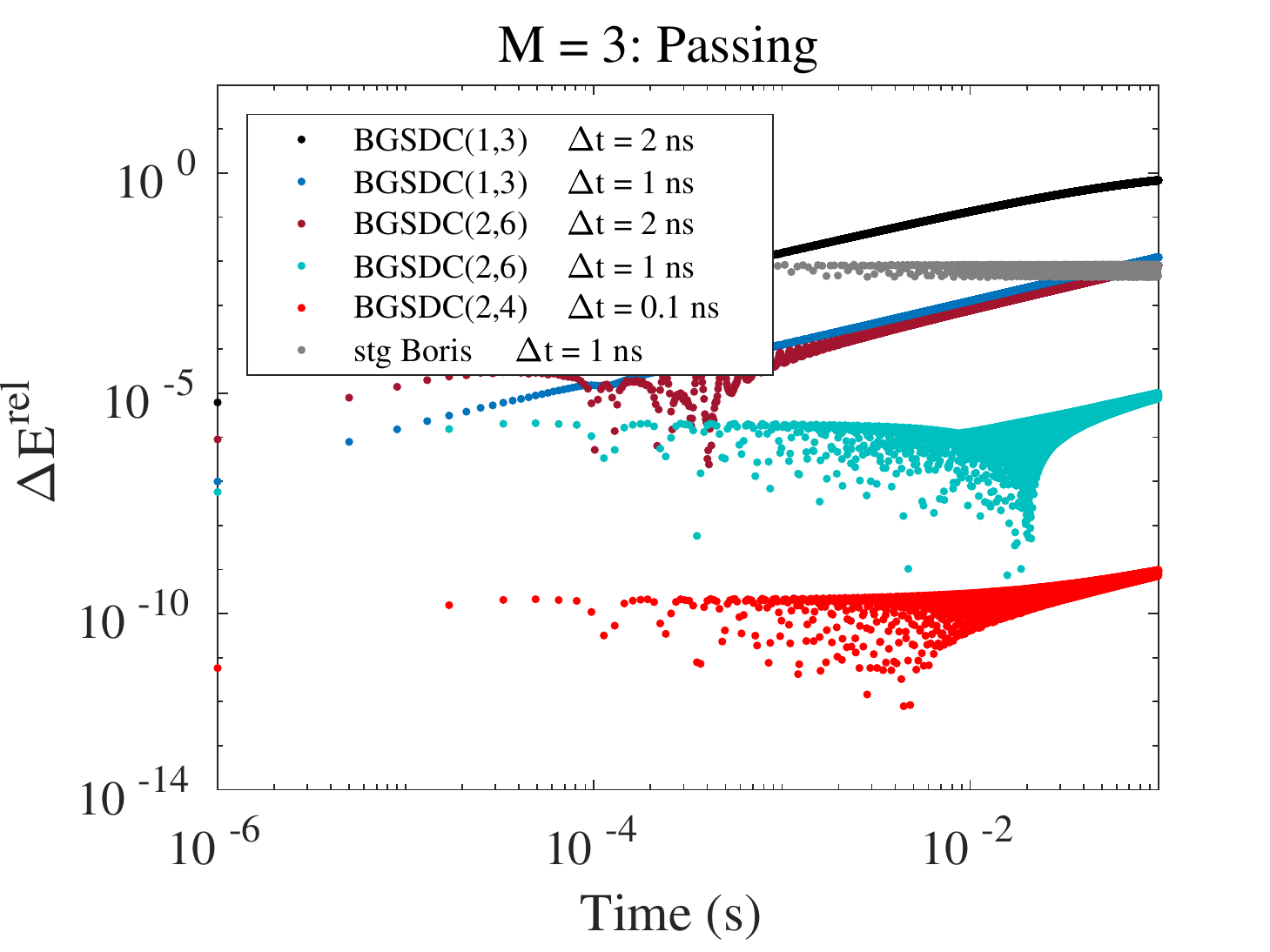} 
	\end{subfigure}
	\hspace{1mm}
	\begin{subfigure}{0.5\textwidth}
		\includegraphics[scale=0.5]{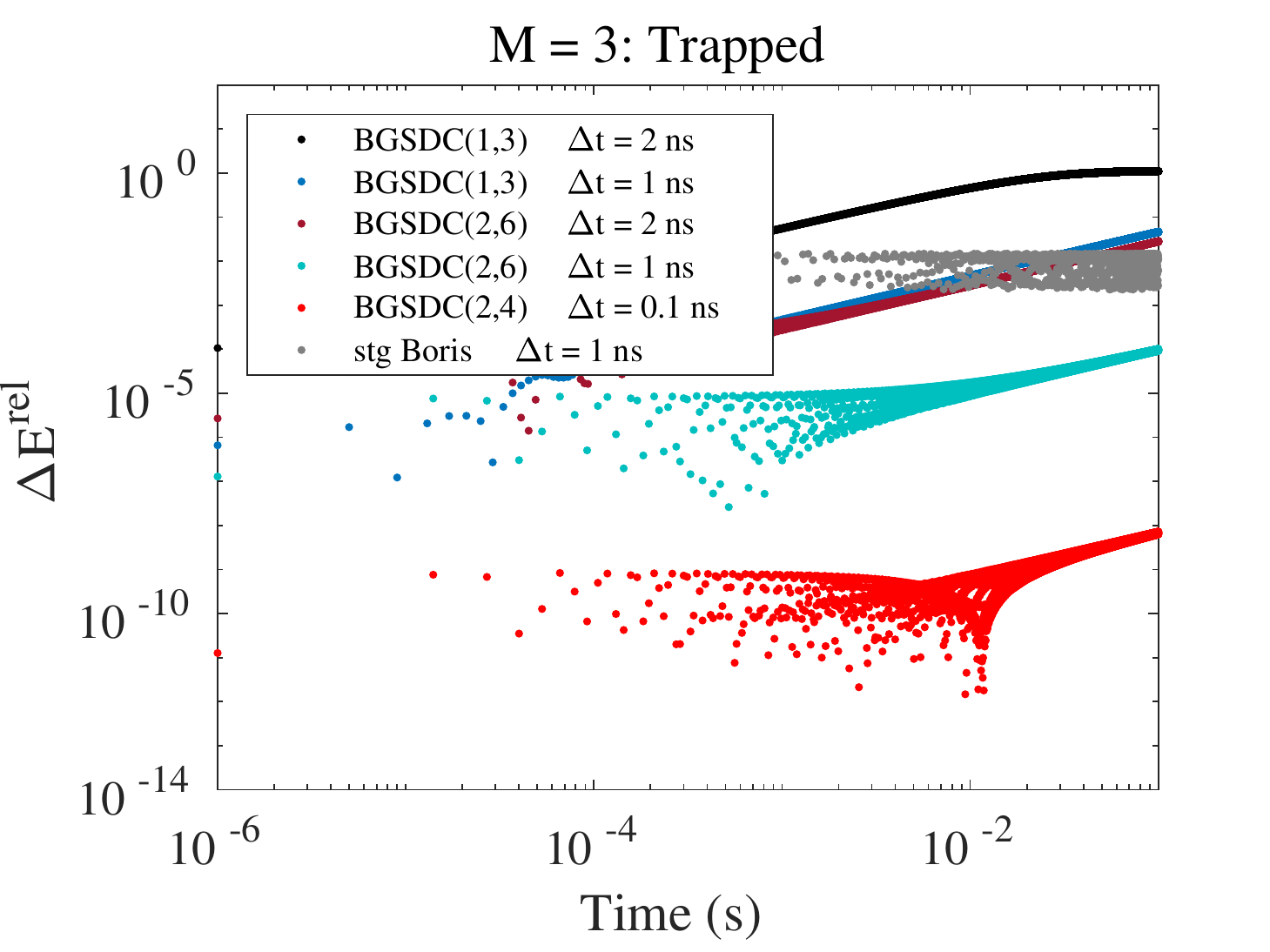}
	\end{subfigure}	
	\begin{subfigure}{0.5\textwidth}
		\includegraphics[scale=0.5]{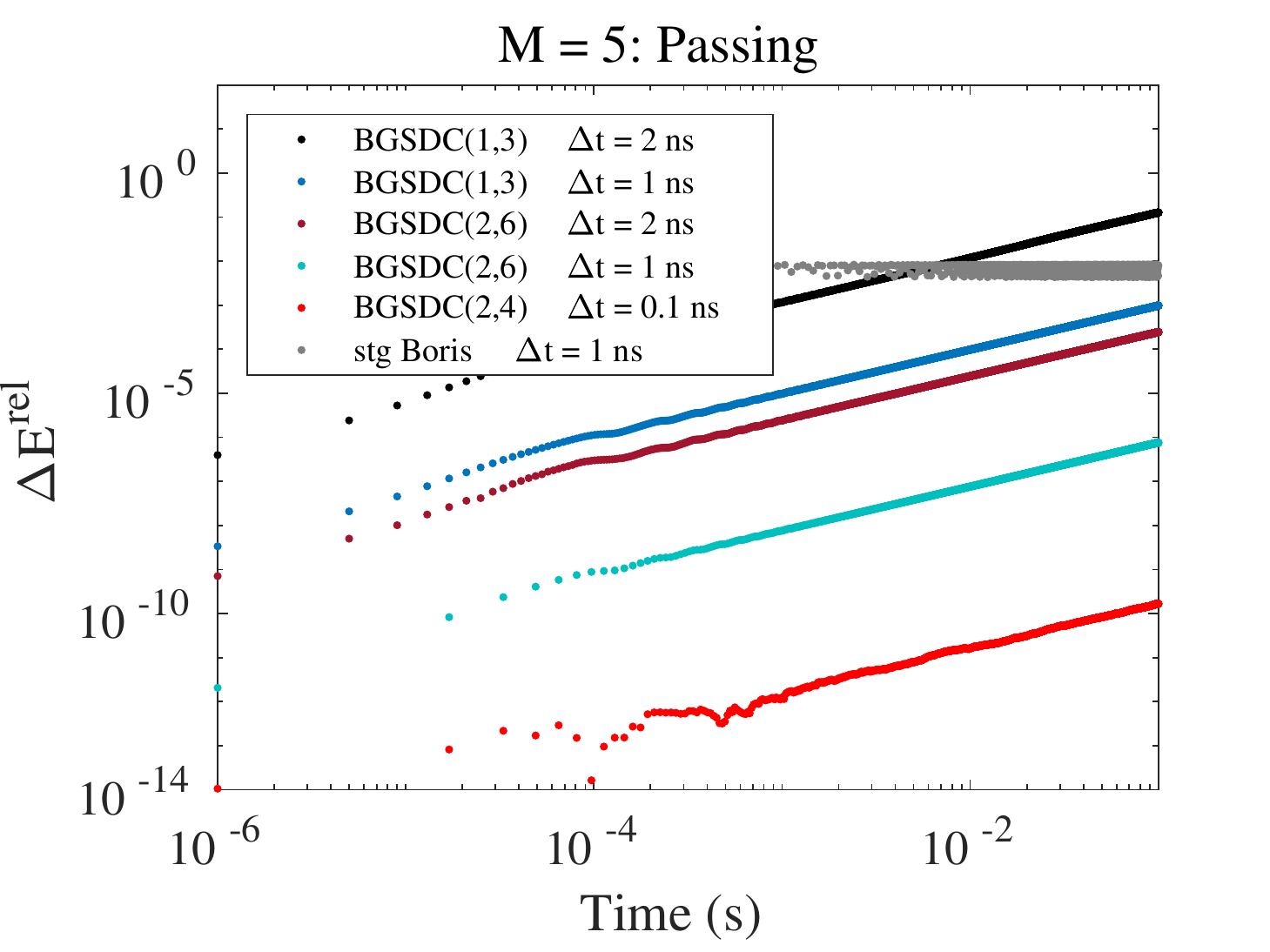} 
	\end{subfigure}
	\hspace{1mm}
	\begin{subfigure}{0.5\textwidth}
		\includegraphics[scale=0.5]{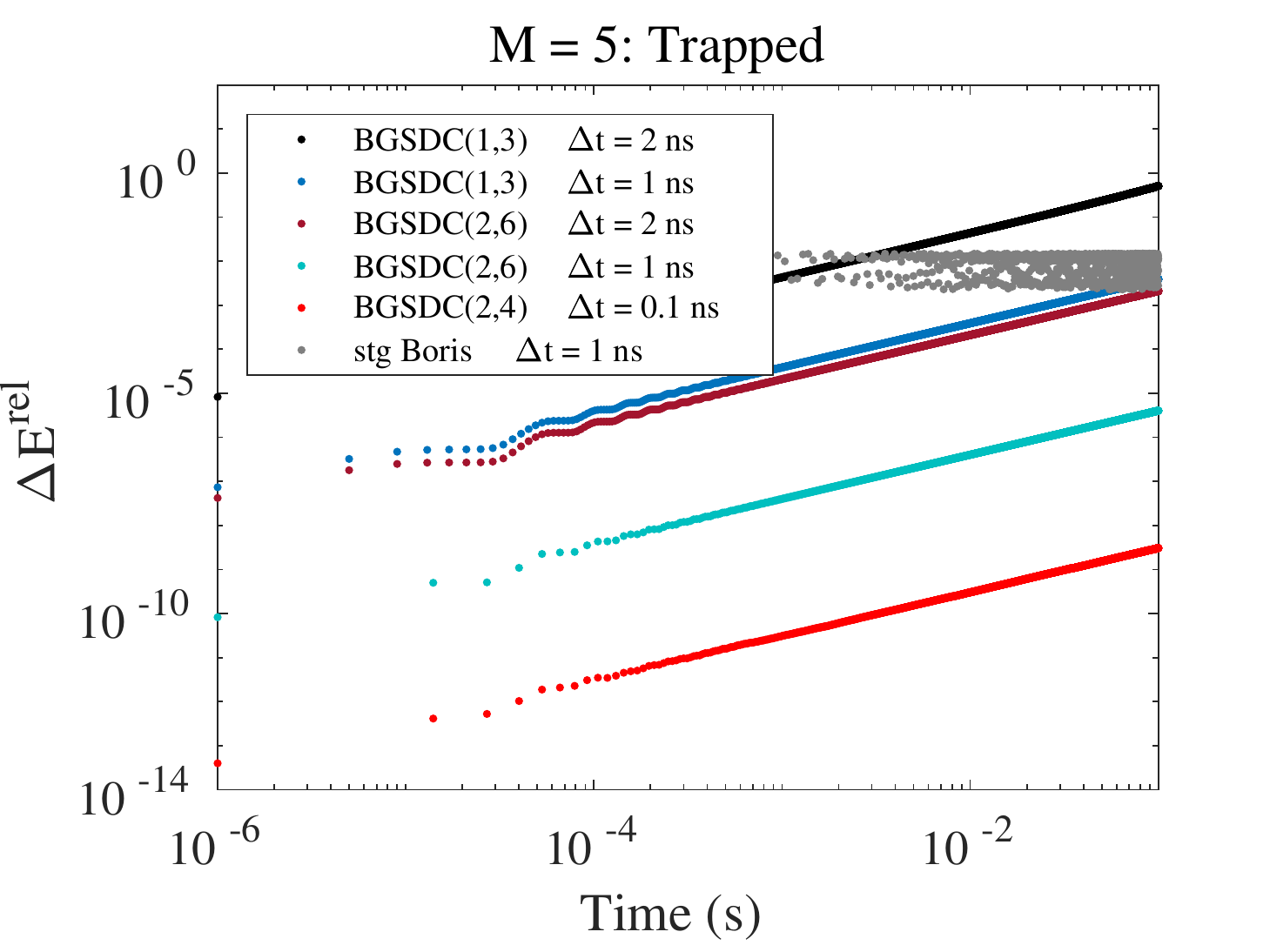}
	\end{subfigure}	
	\caption{Energy error for the passing (left) and trapped (right) particle for different configurations of \revb{BGSDC} for simulations until $t_{\textrm{end}} = \SI{100}{\milli\second}$.}
	\label{fig:solevev_energy_error}
\end{figure}

Figure~\ref{fig:solevev_energy_error} shows the long-time energy error for various configurations of \revb{BGSDC} for the Solev'ev test case.
\revb{Although formally the collocation formulation underlying BGSDC is symmetric (because we use Gauss-Lobatto nodes), accumulation of round-off error still causes energy drift, a well documented problem of methods that rely on iterative solvers~\cite{HairerEtAl2008}.}
However, the growth is relatively mild and energy errors are typically small, even after millions of steps, \revb{if the number of iterations is sufficiently high}.
\revb{Results are similar for the passing and trapped particle.
BGSDC(2,6) with $M=3$ nodes and a time step of \SI{1}{\nano\second} delivers a final error of around $10^{-5}$ for the passing and $10^{-4}$ for the trapped case.
}

%% file: conclusions.tex
\section{Conclusions and future work}
The paper introduces Boris-GMRES-SDC \revb{(BGSDC)}, a new high order algorithm to numerically solve the Lorentz equations based on the widely-used Boris method.
BGSDC relies on a combination of spectral deferred corrections for second order problems and a GMRES-based convergence accelerator originally devised for first order problems.
Since it freezes the magnetic field over the GMRES iterations to linearise the collocation problem, its applicability is limited to cases where the magnetic field does not change substantially over the course of one time step.
Parts of the introduced algorithm are amenable to parallelisation, opening up a possibility to introduce some degree of parallelism in time across the method (\reva{following} the classification by Gear~\cite{Gear1988}), but this is left for future work.

The new algorithm is compared against the standard Boris method for two problems, a magnetic mirror and a Solev'ev equilibrium, the latter resembling the magnetic field of the JET experimental fusion reactor at the Culham Centre for Fusion Energy.
For the Solev'ev equilibrium, our examples show that if precisions in the millimetre range are required, \revb{BGSDC} can reduce computational effort by factors of up to $4$ compared to the standard Boris method.
Gains will be greater for even smaller accuracies and will decrease if less accuracy is needed.
While the break even point from where \revb{BGSDC} cannot produce computational gains is hard to pinpoint, our result\reva{s} suggest it to be for precisions in the centimetre range.
A properly parallelised implementation of \revb{BGSDC} together with an effective adoption of the parameter optimisation strategy by Weiser~\cite{Weiser2014} may still outperform the classical Boris method but is left for future work.